%% file: main.tex
\documentclass{article}
\usepackage{authblk}
\usepackage{blindtext}

\input{packages}

\title{Stable neural networks and connections to continuous dynamical systems}

\author[1]{Matthias J. Ehrhardt}
\author[2]{Davide Murari }
\author[2]{Ferdia Sherry}
\affil[1]{Department of Mathematical Sciences, University of Bath}
\affil[2]{Department of Applied Mathematics and Theoretical Physics, University of Cambridge}

\date{}

\begin{document}
\maketitle
\begin{abstract}
The existence of instabilities, for example in the form of adversarial examples, has given rise to a highly active area of research concerning itself with understanding and enhancing the stability of neural networks. We focus on a popular branch within this area which draws on connections to continuous dynamical systems and optimal control, giving a bird’s eye view of this area. We identify and describe the fundamental concepts that underlie much of the existing work in this area. Following this, we go into more detail on a specific approach to designing stable neural networks, developing the theoretical background and giving a description of how these networks can be implemented. We provide code that implements the approach that can be adapted and extended by the reader. The code further includes a notebook with a fleshed-out toy example on adversarial robustness of image classification that can be run without heavy requirements on the reader’s computer. We finish by discussing this toy example so that the reader can interactively follow along on their computer. This work will be included as a chapter of a book on scientific machine learning, which is currently under revision and aimed at students.
\end{abstract}

It was observed in \cite{szegedy2013intriguing} that tiny, imperceptible perturbations to input images can cause neural networks to misclassify inputs that were previously classified correctly. A remedy to this problem is to make the network stable by controlling the Lipschitz constant of the network \cite{tsuzuku2018lipschitz,sherryDesigningStableNeural2024,trockman2021orthogonalizing,prach20241}. Constraining the Lipschitz constant of neural networks is also fundamental in several data-driven techniques in inverse problems, an area of study that has attracted a lot of interest lately, see, e.g., \cite{benning2018modern,maass2019deep,arridge2019solving}. Along the same line, cleverly designing the network layers of a very deep network is essential for a stable training procedure \cite{galimberti2023hamiltonian,he2015deepresiduallearningimage}. All three examples mentioned in the previous lines have one aspect connecting them: stability. The existence of instabilities in neural networks, such as adversarial examples, has given rise to a highly active area of research focused on understanding and enhancing the stability of neural networks. Building stable neural networks is challenging, since neural networks are highly non-linear parametric functions whose stability properties are hard to understand. To present a common viewpoint to the stability problem, we focus on a popular branch within this area which draws on connections to continuous dynamical systems and optimal control, giving a bird's eye view of this area. We identify and describe the fundamental concepts that underlie much of the existing work in this area. Depending on the application of interest, there are different notions of stability to consider. The goal of this work is to provide an extensive coverage of the most common ones in deep learning: non-expansiveness, Lyapunov stability, and stable network training to avoid vanishing gradient problems. We dedicate a section to each of these notions. We then go into more details on a specific approach to designing stable neural networks with controlled Lipschitz constants, developing the theoretical background and giving a description of how these networks can be implemented.

We provide code that implements this approach, which can be adapted and extended by the reader. The code takes the form of two \texttt{jupyter} notebooks collected in the repository \href{https://github.com/davidemurari/bookChapterDS}{https://github.com/davidemurari/bookChapterDS}. The first is concerned with regularising an ill-conditioned inverse problem, and the second investigates the problem of adversarially robust image classification and the application of the proposed networks for this purpose. The end of the paper describes the problems and methods considered in the code in more detail. We focus on low-dimensional and didactic examples to facilitate visualisation, decrease the time and memory costs of the simulations, and focus on the methodology rather than the application itself. Still, the methods we present extend naturally to higher-dimensional inverse problems and classification tasks, with all stability guarantees preserved. Throughout the manuscript, we also mention some more realistic problems where the presented methodology can be applied.

This manuscript will be included as a chapter in a book on scientific machine learning, currently under revision and aimed at students. The style and the focus on worked-out examples, including exercises, and simple implementations are due to this scope.

\paragraph{Outline of the paper.} This work is structured as follows. Section \ref{sec:intro} provides a detailed motivation of why neural networks suffer from stability problems, and anticipates the solutions that this work expands on. Section \ref{sec:nns_discretised_dynsys} presents the connection between dynamical systems and neural networks that we leverage to formalise the notion of stability and design stable neural networks. The following sections rely on this dynamical systems viewpoint to build networks with specific stability properties. In Section \ref{sec:nonexpansiveness} we show how to build non-expansive networks. In Section \ref{sec:hamiltonian} we present networks that do not suffer from vanishing gradient problems and can hence be trained despite being very deep. This section leverages the formalism of Hamiltonian mechanics to build stable networks. Section \ref{sec:stableEquilibria} studies networks designed to approximate unknown dynamical systems that are known to be Lyapunov stable, and hence presents a third notion of stability and how it relates to deep learning. Section \ref{sec:examples} returns to non-expansive networks and is dedicated to two specific applications in inverse problems and robust image classification. Here, we present the details of a numerical implementation of these networks, reporting and commenting on the results of the discussed simulations. Throughout, we also include some exercises to consolidate the understanding of the content. Still, their resolution is not fundamental to understanding the content.

\input{intro}

\input{nns_discretised_dynsys}
\input{nonexpansiveness}
\input{hamiltonianNetworks}
\input{asymptoticStability}
\input{worked_examples}

\input{acknowledgements}

\bibliographystyle{plain}
\bibliography{references}

\end{document}

%% file: packages.tex
\usepackage{url}
\usepackage{hyperref}
\usepackage{amsmath,amsthm,amsfonts}
\usepackage{mathtools}
\usepackage[export]{adjustbox}
\usepackage{subcaption}

\usepackage[sort,nocompress]{cite}

\newcommand{\new}[1]{#1}

\newcommand{\last}[1]{#1}

\renewcommand{\*}[1]{\mathbf{#1}}
\newcommand{\dd}{\mathrm d}
\newcommand{\R}{\mathbb{R}}
\newcommand{\N}{\mathbb{N}}
\newcommand{\y}{\mathbf{y}}
\newcommand{\x}{\mathbf{x}}
\newcommand{\bu}{\mathbf{u}}

\newcommand{\q}{\mathbf{q}}
\newcommand{\p}{\mathbf{p}}

\DeclareMathOperator{\id}{id}

\DeclareMathOperator*{\argmax}{argmax}

\newcommand{\dx}{\dot{\mathbf{x}}}

\newcommand{\J}{J}

\theoremstyle{definition}
\newtheorem{exercise}{Exercise}
\newtheorem{example}{Example}

\newtheorem{definition}{Definition}

\usepackage{xcolor}
\usepackage{listings}
\usepackage{caption}
\usepackage{newfloat}

\DeclareFloatingEnvironment[fileext=lod,placement={!ht},name=Code Snippet]{listing}
\newenvironment{code}{\captionsetup{type=listing}}{}

\definecolor{pyGreen}{RGB}{0,128,0}     
\definecolor{pyBlue}{RGB}{0,0,180}      
\definecolor{pyGray}{RGB}{102,102,102}  
\definecolor{pyString}{RGB}{163,21,21}  

\lstdefinestyle{pyLikeMinted}{
  language=Python,
  basicstyle=\ttfamily\footnotesize,
  keywordstyle=\bfseries\color{pyGreen},
  commentstyle=\itshape\color{pyGray},
  stringstyle=\color{pyString},
  emphstyle=\bfseries\color{pyBlue},
  showstringspaces=false,
  columns=fullflexible,
  keepspaces=true,
  tabsize=2,
  breaklines=true,
  frame=lines,
  framerule=0.4pt,
  rulecolor=\color{black},
}

\usepackage{csquotes}

%% file: intro.tex
\section{The need for stable neural networks} \label{sec:intro}

Despite the great successes of deep learning in all areas of science and technology, most off-the-shelf neural networks show instabilities: tiny perturbations of the input lead to dramatic consequences in the output. These instabilities may be exploited in adversarial attacks\cite{goodfellow2014explaining,szegedy2013intriguing,madry2017towards} and are particularly problematic in high-risk applications like medical imaging \cite{antun2020instabilities}. In this work, we will study stable neural network architectures designed via the mature mathematical framework of continuous dynamical systems, i.e., differential equations.

Deep neural networks take the form of a function $\Phi : \R^d \to \R^c$,
\begin{align}
    \Phi = \Phi_\ell \circ \dots \circ \Phi_1 ,\label{eq:deepnet}
\end{align}
which comprises of many \textit{layers} $\Phi_i : \R^{n_{i-1}} \to \R^{n_i}, i=1, \dots, \ell$ with $n_0 = d$ and $n_\ell = c$. Specific examples for $\Phi_i$ give particular neural network architectures like the \textit{multi-layer perceptron (MLP)} \cite{rosenblatt1958perceptron} $\Phi_i(\*z) = \sigma(A_i \*z + \*b_i)$ with the \textit{weight matrix} $A_i \in \R^{n_i \times n_{i-1}}$, the bias $\*b_i \in \R^{n_i}$ and the \textit{activation function} $\sigma : \R^{n_i} \to \R^{n_i}$ acts componentwise, e.g., a common choice is the Rectified Linear Unit $[\sigma(\*z)]_i = \max(z_i, 0)$ or the sigmoid $[\sigma(\*z)]_i = 1/(1 + \exp(-z_i))$. 
Since the activation function is always applied componentwise, we do not distinguish between it and the function applied to each component which we also denote by $\sigma:\R\to\R$. If we replace the weight matrix $A_i$ with a convolution, then we speak of a \textit{convolutional neural network (CNN)} \cite{fukushima1980neocognitron, yamashita2018convolutional}.

A common problem in deep neural networks is vanishing and exploding gradients, which prohibit effective training of network parameters. This means that the gradients with respect to the parameters become either very small or very large during training. One strategy proposed in the literature to combat this phenomenon is the ResNet \cite{he2015deepresiduallearningimage} which replaces the individual layers and adds so-called \emph{skip connections}:
\begin{align}
    \Phi_i(\*z) = \*z + \sigma(A_i \*z + \*b_i). \label{eq:resnet}
\end{align}
There is an intrinsic connection between ResNets and the theory of dynamical systems, which we will discuss in more detail in Section \ref{sec:nns_discretised_dynsys}.

Coming back to the topic of stability, the simplest notion of stability is \textit{Lipschitz continuity}. We call a neural network $\Phi$ \textit{stable} if it is $L$-Lipschitz continuous, i.e., there exists a constant $L \geq 0$ such that
\begin{align}
    \|\Phi(\*x) - \Phi(\*y)\|_2 \leq L \|\*x - \*y\|_2\new{,\,\,\forall \x,\y\in\R^d}. \label{eq:Lipschitz}
\end{align}
Due to the layered structure of deep neural networks \eqref{eq:deepnet}, we can relate the Lipschitz constant $L$ to the Lipschitz constants of the individual layers $L_i$ as $L \leq \Pi_{i=1}^\ell L_i$ \cite{gouk2021regularisation}. It was argued in \cite{bungert2021clip} that such an estimate is pessimistic and hinders practical usefulness.

Stability of neural networks is desirable in many contexts such as the stable solution to inverse problems, classification that is robust to errors, efficient training of deep neural networks. It is also used in the context of generative models and is frequently used for (Wasserstein) GANs to regularise the discriminator~\cite{pmlr-v70-arjovsky17a}, for instance via spectral normalisation \cite{miyatoSpectralNormalizationGenerative2018}.

We now expand on three specific examples which describe potential use cases of stable neural networks.

\begin{example}[Inverse Problems] The first example we consider is inverse problems where we are interested in recovering some quantity $\x^\dagger \in \R^d$ from measurements $\y^\delta = A\x^\dagger + \*z \in \R^m$ where $\*z$ is some measurement noise.
There are numerous applications where such modelling is useful, such as X-ray computerised tomography in medical imaging or material science. Since the measurements $\y^\delta$ contain noise and the matrix $A$ is usually ill-conditioned, simply \enquote{inverting} $A$ does not lead to useful solutions. This is illustrated in Figure \ref{fig:ip}, where we consider the simple yet insightful example with
$$ A = \begin{pmatrix} 1 + \varepsilon & 1 \\ 1 & 1 + \varepsilon\end{pmatrix}$$
for $\varepsilon = 1/4$. The figure shows that an inversion of $A$ is only meaningful in the absence of noise.
Its eigenvalues are given by $2 + \varepsilon$ and $\varepsilon$, thus its condition number $1 + 2/\varepsilon$ making the problem severely ill-conditioned for small $\varepsilon$. Similar to classical regularisation theory, the inverse of $A$ can be approximated with a stable neural network, $\Phi \approx A^{-1}$, making the reconstruction reliable even for noisy measurements. This will be considered in more detail in Section \ref{sec:inverse_problem}. \new{Our discussion of this example is self-contained. For further material on inverse problems and related data-driven techniques, see \cite{benning2018modern,arridge2019solving}.}

\begin{figure}[!ht]
    \centering
    \includegraphics[scale=1]{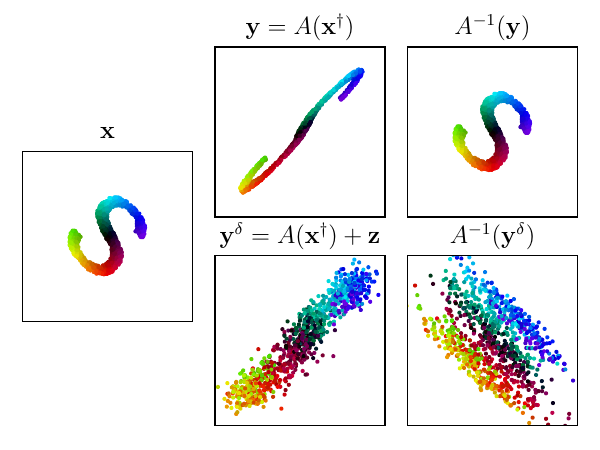}
    \caption{An illustration of how the presence of noise and ill-conditioning of the matrix $A$ combine to complicate the inversion process.} \label{fig:ip}
\end{figure}
\end{example}

\begin{example}[Classification]
The second example we consider is the problem of classifying inputs $\new{\x} \in \R^d$ into $c \in \mathbb N$ classes. A standard approach to this problem models the classifier using a neural network $\Phi: \R^d \to \R^c$ and predicts the class for an input $\*x$ as $\argmax_{k\in \{1,\ldots, c\}} \Phi_k(\*x)$. The outputs of $\Phi$ may be interpreted as \emph{logits}, meaning that $\exp(\Phi_k(\*x)) / \sum_{i=1}^c \exp(\Phi_i(\*x))$ is treated as a probability that $\*x$ is of class $k$. Regardless of the interpretation, Lipschitz continuity of $\Phi$ can be used to \emph{certify the robustness} of the predictions. 

Associated with $\Phi$, we can define the predicted class $\hat k : X \to \{1, \ldots, c\}$ as $\hat k(\new{\x}) = \argmax_{k\in \{1,\ldots, c\}} \Phi_k(\*x)$ and the \emph{margin} $m : \R^d \to \R$ as 
\[m(\*x) = \Phi_{\hat k(\*x)}(\*x)  - \max_{k \in \{1,\ldots, c\} \setminus \{\hat k(\*x)\}} \Phi_k(\*x).\]
One can think of the margin as some wiggle room in the accuracy of the prediction. Even if the predicted value shrinks to the margin, the model's prediction remains the same. If the classifier is stable, then there can be errors in the data that do not alter the classification result. In detail, if $\Phi$ is $L$-Lipschitz, then the classification is robust, in the sense that any $\*y \in \R^d$ with 
$$\|\*x - \*y\|_{\new{2}}< \frac{m(\*x)}{2L}$$ 
is predicted to be of the same class as $\*x$, i.e., $\hat k(\*x) = \hat k(\*y)$. This result can be used to give robustness guarantees for existing classifiers, but can also be used to motivate the training of robust classifiers: if we can upper bound the Lipschitz constant of a neural network.

In order to justify the statement above, one can show that 
\[\Phi_{\hat k(\*x)}(\*y) - \max_{k \in \{1, \ldots, c\}\setminus \{\new{\hat{k}}(\*x)\}} \Phi_k(\*y) > \Phi_{\hat k}(\*x) - \frac{m(\*x)}{2} + \frac{m(\*x)}{2} - \Phi_{\hat k}(\*x) = 0,\]
showing that $\hat k(\*y):=\argmax_{k\in \{1,\ldots, c\}} \Phi_k(\*y) = \hat k(\*x)$, so that the predictions of $\Phi$ at $\*x$ and $\*y$ match. 

\begin{exercise}
Fill in the details of the argument above. In particular, show that 
\[\Phi_{\hat k(\*x)}(\*y) + \frac{m(\*x)}{2} > \Phi_{\hat k(\*x)}(\*x),\]
and
\[\max_{k \in \{1,\ldots, c\} \setminus \{\hat k(\*x)\}} \Phi_{k}(\*x) + \frac{m(\*x)}{2} > \max_{k \in \{1,\ldots, c\} \setminus \{\hat k(\*x)\}} \Phi_{k}(\*y).\]
\end{exercise}

We will come back to image classification in Section \ref{sec:robust_image_classification}. 
\end{example}

\begin{example}[Stable Training] \label{ex:training}
The third example we want to consider in more detail is the stability of the forward propagation in the context of the network training procedure. For supervised learning, the network may be trained by minimising a loss function
\[
    \mathcal{L}\left(\theta\right)=\frac{1}{n}\sum_{i=1}^n\left\|\Phi_{\theta}\left(\x^i\right)-\y^i\right\|_2^2,
\]
given a dataset $\left\{\left(\x^i,\y^i\right)\right\}_{i=1}^n \subset \R^d \times \R^c$. Here we make the dependency of the network $\Phi$ on its parameters $\theta$ explicit by writing $\Phi_\theta$. The process of training the neural network $\Phi_\theta$ involves the computation of the gradients of the loss function $\mathcal{L}$ with respect to the network weights $\theta$. For example, if the parameters are trained by gradient descent, then the $i$th component of the parameter vector of the $j$th layer $\theta_j$, which we denote as $\theta_{ij}$, is updated as
\begin{equation}\label{eq:gradDesc}
    \theta_{ij}^{k+1}=\theta_{ij}^k - \tau^k \partial_{\theta_{ij}}\mathcal{L}\left(\theta\right) = \theta_{ij}^k - \frac{\tau^k}{n}\sum_{m=1}^n\partial_{\theta_{ij}}\mathcal L_m\left(\theta\right),
\end{equation}
where we used the notation $\mathcal L_m\left(\theta\right)=\left\|\Phi_{\theta}\left(\x^m\right)-\y^m\right\|_2^2$, and $\tau^k$ is the step-size, also called the learning rate in this context. 

As seen earlier, a neural network with $\ell$ layers is given by $\Phi_{\theta} = \Phi_{\theta_\ell}\circ \cdots \circ \Phi_{\theta_1}$. Alternatively, we can write it as $\Phi_\theta(\*x) = \*x_{\ell+1}$ with
\begin{align*}
    \*x_{t+1} = \Phi_{\theta_t}(\*x_t), \quad t = 1, 2, \dots, \ell,
\end{align*}
and $\*x_{1} = \*x$. One may notice that, especially for many layers $\ell$, the compositional nature of $\Phi_{\theta}$ can lead to vanishing gradients. Indeed, by the chain rule, we see that for any $m$
\begin{equation*}
    \partial_{\theta_{ij}}\mathcal L_m = \langle \partial_{\x_{j+1}^m} \mathcal L_m,  \partial_{\theta_{ij}} \x_{j+1}^m\rangle = \Big\langle\Big(\prod_{t=j+1}^\ell \partial_{\x_t^m}\x_{t+1}^m \Big)\partial_{\x_{\ell+1}^m} \mathcal L_m, \partial_{\theta_{ij}} \x_{j+1}^m\Big\rangle,
\end{equation*}
where $\langle\*x, \*y\rangle$ is the Euclidean inner product between two vectors $\*x, \*y$. Together with the inequality
\begin{equation}\label{eq:normIneq}
    \Bigg\|\prod_{t=j+1}^\ell \partial_{\x_{t}^m} \x_{t+1}^m\Bigg\|_2\leq \prod_{t=j+1}^\ell\Big\|\partial_{\x_{t}^m} \x_{t+1}^m \Big\|_2,
\end{equation}
imply that if $\ell$ is large and the norms on the right of \eqref{eq:normIneq} are smaller than $1$, then the gradient $\nabla_{\theta_{ij}}\mathcal L_m$ will be very small (or converge to zero for $\ell \to \infty$), hence leading to the impossibility of updating the weights in a meaningful way using gradient information as in \eqref{eq:gradDesc}. This is illustrated in Figure \ref{fig:ex:training}, where the vanishing gradient phenomenon leads to poor classification results.
\begin{figure}[ht!]
\centering
\includegraphics[scale=1]{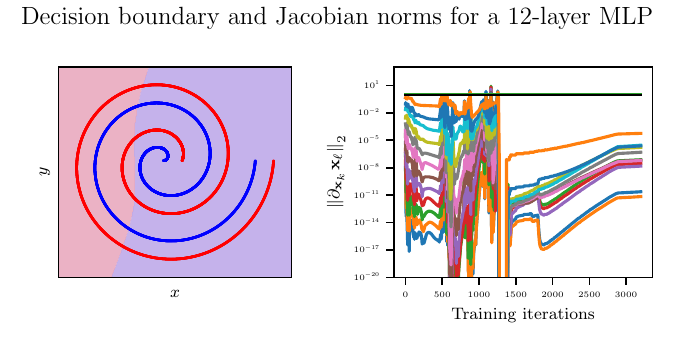}
\caption{MLP network with 12 layers trained to distinguish red from blue points. On the left, the learned decision boundary cannot accurately separate the points, resulting in a test accuracy of 51\%. On the right, the norms of the Jacobians through the training iterations for a fixed data point, showing a severe attenuation of information as we progress through the network, an issue known as the vanishing gradient problem.
}
\label{fig:ex:training}
\end{figure}

Because of this fundamental issue, it is important to suitably design the layers $\Phi_{\theta_1},\ldots,\Phi_{\theta_\ell}$, so that
$\|\partial_{\x^m_j} \x_{j+1}^m\|_2
$ is of moderate size. We will revisit this problem in Section \ref{sec:hamiltonian}.
\end{example}

%% file: nns_discretised_dynsys.tex
\section{{Neural networks as discretised dynamical systems}}
\label{sec:nns_discretised_dynsys}
The development of calculus by Newton and Leibniz in the 17\textsuperscript{th} century went hand in hand with its applications in the mathematical modelling of mechanical systems. Subsequently, various interconnected subfields have been developed, including Lagrangian and Hamiltonian mechanics, started by Lagrange in the 18\textsuperscript{th} century and Hamilton in the 19\textsuperscript{th} century, respectively, and the systematic and far-reaching study of such and other models in the theory of dynamical systems, initiated by Poincar\'e in the 19\textsuperscript{th} and 20\textsuperscript{th} centuries. The study of dynamical systems has become indispensable for modern science and engineering, and we will study how various ideas from this field of research can be used to great effect to aid in the design and understanding of neural networks.

The dynamical systems that we will focus on are described by ordinary differential equations (ODEs), although we note that extensions are possible to other classes of models, most notably including partial differential equations (PDEs), corresponding to infinite-dimensional state spaces, and stochastic differential equations (SDEs), which are driven by stochastic processes in addition to deterministic forces.

We are particularly interested in so-called \emph{initial value problems} (IVPs). Given a differential equation, embodied in its vector field $X:[0,T]\times \R^d\to\R^d$, an initial time $t_0\in \R$ and an initial condition $\x_0 \in \R^d$, the goal is to determine a trajectory $\x:[t_0,T] \to \R^d$ through $\x_0$:
\begin{equation}
    \begin{cases}
        \dot\x(t) = X(t, \x(t)), \\
        \x(t_0) = \x_0.
    \end{cases}
    \label{prob:ivp}
\end{equation}
Here, and in what follows, we will use dot notation to indicate derivatives with respect to a time variable, e.g., $\dot \x = \dd \x/\dd t$. \new{We now present a few basic results on IVPs. To further study this topic, we refer to \cite{doi:10.1137/1.9781611975161,strogatz2018nonlinear}.}
\begin{exercise}
    \label{exercise:higher_order_ode}
    The form of the differential equation in \eqref{prob:ivp} may seem to be overly restrictive, as it contains only first-order derivatives. As an example, Newton's equations of motion (in their standard form) are second-order ODEs: $\ddot\x(t) = - \nabla V(\x(t))/m$. Show that, by appropriately augmenting the state space, we can reinterpret a higher-order ODE,
    \begin{equation}
        \x^{(n)}(t) = f(t, \x(t), \dot\x(t), \ldots, \x^{(n - 1)}(t)),
    \end{equation}
    for some $n\in \N$, in the form of \eqref{prob:ivp}. Here, we denote by $\x^{(n)}$ the $n$-th derivative of $\x$.
\end{exercise}

Under standard assumptions on the vector field $X$, we can guarantee the (local) existence and uniqueness of solutions to \eqref{prob:ivp}: if $X$ is locally Lipschitz-continuous \new{in the second argument}, the \emph{Picard--Lindel\"of theorem} gives such a guarantee. A large range of types of dynamics, each with their own characteristic behaviours, can be described in the form of \eqref{prob:ivp}. For example, suppose $X$ is the negative gradient of a convex potential. In that case, we get \emph{non-expansive} (i.e., 1-Lipschitz continuous) dynamics as will be described in Section \ref{sec:nonexpansiveness}, while if $X$ is a Hamiltonian vector field, the resulting dynamics conserves energy, as will be discussed in Section \ref{sec:hamiltonian}.

From now on, the question of the existence and uniqueness of solutions will be of no concern, as the models that we are considering are well-behaved in this respect. For convenience, then, we will assume that $t_0 = 0$, and denote the solution to \eqref{prob:ivp}, evaluated at a time $t$, by $\phi^{t}_X(\x_0)$. In particular, if the vector field does not depend on time, in which case it is said to be \emph{autonomous}, this gives us a continuous group of transformations: $\phi_X^0 = \id$ and $\phi_X^{t + s} = \phi_X^t\circ \phi_X^s$. \new{For the simple case $X(t,\x)=A\x$, $A\in\R^{d\times d}$, one obtains $\phi^t_X(\x_0)=\exp(At)\x_0$, for any $t\geq 0$ and $\x_0\in\R^d$, where $\exp(A)=\sum_{k\geq 0}A^k/k!$ stands for the matrix exponential.}

\subsection{Discretising ordinary differential equations}
For all but the simplest ODEs, it is impossible to explicitly solve problem~\eqref{prob:ivp}. As a result, it becomes essential to numerically approximate its solution, a topic which received some attention in the early days of calculus but which has exploded in interest in the past century with the advent of computers and the associated scaling up of scientific problems to be tackled.

\begin{figure}[!htb]
    \centering
    \includegraphics[scale=1]{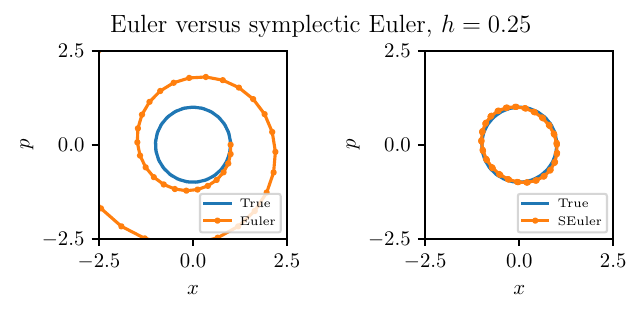}
    \caption{A demonstration of the behaviours discussed in Exercise \ref{exercise:ho_euler_symplectic} for the harmonic oscillator, with an initial condition $x_0 = 1, p_0=0$.}
    \label{fig:ho_euler_symplectic}
\end{figure}

Many such methods can be classified as time-stepping methods, which approximate the solution trajectory of \eqref{prob:ivp} at discrete points in time by sequentially composing approximations to the true time steps. The simplest example of such a method is the \emph{explicit Euler} method, also known as the forward Euler method, or simply the Euler method: we proceed by taking a first-order Taylor expansion of the solution,
\begin{equation}
    \x(t + h) = \x (t) + h\dot \x (t)   + O(h^2) =  \x(t) + hX(t, \x(t)) + O(h^2),
    \label{eq:taylor_euler_ode}
\end{equation}
and simply drop the higher-order terms. Hence, given an initial time $t_0$, we can approximate the solution at discrete times, $\{\x(t_0 + n h)\}_{n=0}^N$, by the sequence recursively defined as follows:
\begin{equation}
    \begin{cases}
        \x^\textrm{Euler}_0        := \x_0 , \\
        \x^\textrm{Euler}_{n + 1}  := \x^\textrm{Euler}_{n} + h X(t_0 + nh, \x^\textrm{Euler}_n).
    \end{cases}
    \label{eq:euler_method}
\end{equation}
The Taylor expansion \eqref{eq:taylor_euler_ode} shows that this method is \emph{consistent}, in the sense that the local error made by a single step of the method vanishes as $h\to 0$. It is a fundamental theorem of numerical analysis that\new{, for a stable method, consistency} implies a global \emph{convergence} result too: with a fixed time horizon, the global error incurred by the Euler method in approximating the true trajectory is of order $h$ as $h\to 0$. \new{ More details on numerical methods for ODEs can be found in \cite{hairer1993solving}.} Consistency and convergence of a numerical method for ODEs could be considered necessary conditions for its admissibility, but these conditions are far from sufficient to guarantee that the method will perform well in the non-asymptotic setting, where the step size can not be taken to 0, as highlighted in Exercise \ref{exercise:ho_euler_symplectic} and Figure \ref{fig:ho_euler_symplectic}. In this setting, one can consider the use of \emph{structure-preserving numerical methods}, which similarly approximate the solution to \eqref{prob:ivp}, but do so while preserving some of its structural characteristics, such as \emph{symplecticity} (as we will see in Exercise \ref{exercise:ho_euler_symplectic} and Section \ref{sec:hamiltonian}), \emph{conserved quantities}, \emph{dissipation} and \emph{non-expansiveness} (as we will focus on in Sections \ref{sec:nonexpansiveness} and \ref{sec:examples}).

\begin{exercise}[On the importance of structure-preserving numerical methods]
    \label{exercise:ho_euler_symplectic}
    Consider the simple harmonic oscillator, with ODE
    \begin{equation*}
        \begin{pmatrix} \dot x(t)\\\dot p(t)\end{pmatrix} = \begin{pmatrix} p(t)\\-x(t)\end{pmatrix}.
    \end{equation*}
    For this system, we can define a total energy $E(x, p):= (x^2 + p^2)/2$. Assume that $(x_0, p_0)$ is an arbitrary initial condition and consider the associated IVP \eqref{prob:ivp}.
    \begin{itemize}
        \item Show that the energy of the true trajectory is conserved: $E(x(t), p(t)) = E(x_0, p_0)=: E_0$ for all $t \geq 0$.
        \item Show that, when we approximate the trajectory by Euler steps with a step size $h >0$, as in \eqref{eq:euler_method}, the energy behaves as follows, for any $n\in\N$,
              \begin{equation*}
                  E_n:= E(x^{\textrm{Euler}}_n, p^{\textrm{Euler}}_n) = E_0(1 + h^2)^n .
              \end{equation*}
              In particular, we have $E_n \sim E_0 \exp(h^2n)$ as $n\to\infty$, i.e., the energy of the approximate trajectory diverges exponentially.
        \item Modify the Euler integrator \eqref{eq:euler_method} as follows, keeping the initialisation as is:
              \begin{equation}
                  \begin{cases}
                      x^{\textrm{SEuler}}_{n + 1}:= x^{\textrm{SEuler}}_n + h p^{\textrm{SEuler}}_n, \\
                      p^{\textrm{SEuler}}_{n+1}: = p^{\textrm{SEuler}}_n - h x^{\textrm{SEuler}}_{n+1}.
                  \end{cases}
                  \label{eq:ho_symplectic_euler}
              \end{equation}
              Show that with this choice of integrator (which has the same order of approximation as the Euler method) and a step size $0 < h < 2$, the energy along the trajectory, $E_n:= E(x_{\mathrm{SEuler}}(nh), p_{\mathrm{SEuler}}(nh))$, is bounded above and below by small perturbations of the true energy $E_0$. This integrator is known as the symplectic Euler method and will be considered in more detail in Section \ref{sec:hamiltonian}. \textbf{Hint:} recognise that \eqref{eq:ho_symplectic_euler} is a linear update and study the eigenvalues of the corresponding matrix.

    \end{itemize}
\end{exercise}

\last{We now show how these numerical methods can be used to design neural network layers.}

\subsection{\last{From numerical methods to neural networks}\label{se:numericalNetworks}}
Looking at \eqref{eq:euler_method}, a natural connection between neural networks and ODEs arises: the Euler integrator approximates a trajectory by composing simple updates, each of which takes the form of the identity plus a small perturbation. The Euler step takes essentially the same form as a ResNet layer: ResNets replace $hX$ by a (simple) neural network, the weights of which are learned using gradient-based optimisation. Indeed, recall from \eqref{eq:resnet} that in its most basic form a layer of the ResNet is given by
\[\Phi_i(\mathbf x) = \mathbf x + \sigma (A_i\mathbf x + \mathbf b_i),\]
which is to be compared with the update of the Euler integrator in \eqref{eq:euler_method}. ResNets were not initially designed with this connection in mind, but rather with the intention of enabling the training of very deep neural networks by mitigating the vanishing gradient problem, see Example \ref{ex:training}. Having established this connection, however, there are many directions in which the design can be refined to better suit specific problems, including by changing the structure of the parametrised vector fields and by changing the numerical integrator to an integrator that better respects the structure of the system under consideration.

It is worth noting that the building blocks of ODE-based neural networks naturally map between a space $\R^n$ and itself, meaning in particular that the dimensionality of the inputs and outputs must be the same. This may seem overly restrictive: for instance, in the case of image classification, where the goal is to reduce a potentially high-resolution image into a vector, it is necessary to reduce the dimensionality of the intermediate states as they progress through the network. In this setting, it is also common to blow up the number of channels, actually increasing dimensionality, in the first layer. Let us remark that such behaviour can still be obtained using ODE-based neural networks, by interspersing the basic blocks with simple (linear, for example) \emph{lifting} layers, if an increase in dimensionality is needed, or \emph{projection} layers, if a decrease in dimensionality is needed.

\new{Although this work centres on ResNets, it is worth briefly introducing a related architecture that draws directly on ODEs: \textit{neural ordinary differential equations }(Neural ODEs)~\cite{chen2018neural,kidger2022neural}. 
A Neural ODE typically corresponds to the flow map up to a chosen final time, typically $T = 1$, of an ODE parametrised by a neural network.  Thanks to this continuous-time viewpoint, Neural ODEs have become a popular backbone for modern generative-modelling methods.} 

%% file: nonexpansiveness.tex
\section{{Non-expansive neural networks}}
\label{sec:nonexpansiveness}
As mentioned in Section \ref{sec:intro}, Lipschitz continuity is a standard way to quantify the stability of a function. The notion of non-expansiveness extends also to dynamical systems, where we say that a vector field $X:\R^d\to\R^d$ is non-expansive if its flow map $\phi_X^t:\R^d\to\R^d$ is non-expansive for every time $t\geq 0$, i.e., $\|\phi_X^t(\x)-\phi_X^t(\y)\|_2\leq \|\x-\y\|_2$.

Since the flow map is not usually accessible, this definition is not so practical. However, supposing $X$ is sufficiently smooth, we can get a much more practical characterisation of non-expansive dynamical systems by Taylor expansion. Indeed, let us consider a small enough scalar $h$ and consider
\begin{align*}
\phi^{t+h}_X(\x)&=\phi^t_X(\x) + hX(\phi^t_X(\x))+\mathcal{O}(h^2),\\
\phi^{t+h}_X(\y)&=\phi^t_X(\y) + hX(\phi^t_X(\y))+\mathcal{O}(h^2),
\end{align*}
for an arbitrary pair of points $\x,\y\in\R^d$. Then, we see that
\begin{align*}
\left\|\phi^{t+h}_X(\y)-\phi^{t+h}_X(\x)\right\|_2^2 &= \left\|\phi^t_X(\y)-\phi^t_X(\x)\right\|_2^2\\
&+ 2h\langle X(\phi^t_X(\y))-X(\phi^t_X(\x)),\phi^t_X(\y)-\phi^t_X(\x)\rangle + \mathcal{O}(h^2),
\end{align*}
and hence
\begin{align*}
\frac{\dd}{\dd t}\left\|\phi^t_X(\y)-\phi^t_X(\x)\right\|_2^2
&= \lim_{h\to 0}\frac{\left\|\phi^{t+h}_X(\y)-\phi^{t+h}_X(\x)\right\|_2^2-\left\|\phi^t_X(\y)-\phi^t_X(\x)\right\|_2^2}{h}\\
&= 2\langle X(\phi^t_X(\y))-X(\phi^t_X(\x)),\phi^t_X(\y)-\phi^t_X(\x)\rangle.
\end{align*}
This derivation implies that if for every $\x,\y\in\R^d$ one has
\begin{equation}\label{eq:osl}
 \langle X(\y)-X(\x),\y-\x\rangle \leq \nu\|\y-\x\|_2^2,
\end{equation}
then it follows
\begin{align}\label{eq:beforeGronwall}
&\frac{\dd}{\dd t}\left\|\phi^t_X(\y)-\phi^t_X(\x)\right\|_2^2\leq 2 \nu\|\phi^t_X(\y)-\phi^t_X(\x)\|_2^2.
\end{align}
If we define $g(t):=\left\|\phi^t_X(\y)-\phi^t_X(\x)\right\|_2^2$ and multiply both sides of \eqref{eq:beforeGronwall} by the positive scalar $e^{-2\nu t}$, we see that
\[
\frac{\dd}{\dd t}\left(e^{-2\nu t}g(t)\right) = e^{-2\nu t}\dot{g}(t)-2\nu e^{-2\nu t}g(t) \leq 0.
\]
We can thus conclude that $e^{-2\nu t}g(t)$ is monotonically non-increasing, so that $e^{-2\nu t}g(t)\leq g(0)$, and hence we have
\begin{equation}\label{eq:boundTraj}
\left\|\phi^t_X(\y)-\phi^t_X(\x)\right\|_2 \leq e^{\nu t}\left\|\y-\x\right\|_2
\end{equation}
for every $t\geq 0$, $\x,\y\in\R^d$. We remark that the distance between any pair $\x$ and $\y$ is not expanded by the flow map $\phi^t_X$ for $t\geq 0$ whenever $\nu\leq 0$. This analysis motivates the introduction of the following definition. 
\begin{definition}[One-sided Lipschitz inequality]
The vector field $X:\R^d\to\R^d$ is one-sided Lipschitz continuous if it satisfies \eqref{eq:osl} for a scalar $\nu\in\R$ and any pair $\x,\y\in\R^d$. $X$ is a non-expansive vector field if \eqref{eq:osl} holds for a $\nu\leq 0$. \new{$X$ is a contractive vector field if \eqref{eq:osl} holds for a $\nu< 0$}.
\end{definition}

Before moving on, we remark that contractivity can be a pretty restrictive assumption on the dynamics. For example, one can see that a contractive dynamical system has to have a unique asymptotically stable equilibrium point \new{(see also Section~\ref{sec:stableEquilibria} for more details about this concept)}. To verify this behaviour, let $\phi^1_X:\R^d\to\R^d$ be the time-$1$ flow of the contractive vector field $X:\R^d\to\R^d$. Banach's fixed point theorem guarantees that $\phi^1_X$ admits a unique fixed point $\x^*\in\R^d$ such that $\phi^1_X(\x^*)=\x^*$. In case this is an equilibrium point of $X$, it has to be asymptotically stable since for any $\x\in\R^d$ we have
\[
\lim_{t\to+\infty}\left\|\phi^t_X(\x) - \x^*\right\|_2 = \lim_{t\to+\infty}\left\|\phi^t_X(\x) - \phi^t_X(\x^*)\right\|_2  \leq \lim_{t\to+\infty}e^{\nu t}\|\x-\x^*\|_2= 0.
\]
If $\x^*$ is not an equilibrium point, then it has to be part of a periodic orbit of period $1$. This is impossible since the existence of such a periodic orbit would lead to infinitely many fixed points for $\phi^1_X$, allowing us to conclude that, in fact, $\x^*$ must be an equilibrium point.

Even though the condition in \eqref{eq:osl} is more practical than where we started from, it can sometimes be hard to verify. For continuously differentiable vector fields, one can simplify the condition to an equivalent characterisation based on the Jacobian matrix $\partial_{\x} X(\x)\in\R^{d\times d}$. Indeed, by the mean value theorem, for every $\x,\y\in \R^d$, there is $\mathbf{z}=s\x+(1-s)\y$, for some $s\in (0,1)$, such that
\[
X(\y)-X(\x)=\partial_{\x}X(\mathbf{z})(\y-\x).
\]
Thus, \eqref{eq:osl} can be formulated as an equivalent condition
\[
\sup_{\x\in\R^d, \mathbf{v}\in\R^d\setminus\{\boldsymbol{0}\}}\frac{\langle \partial_{\x}X(\x)\mathbf{v},\mathbf{v}\rangle }{\|\mathbf{v}\|_2^2}\leq \nu,
\]
or, equivalently, as
\begin{equation}\label{eq:lmax}
\sup_{\x\in\R^d}\lambda_{\max}\left(\frac{\partial_{\x}X(\x)^\top+\partial_{\x}X(\x)}{2}\right)\leq \nu,
\end{equation}
where $\lambda_{\max}(A)$ is the \new{largest} eigenvalue of some matrix $A$.
\begin{exercise}
This exercise relates the one-sided Lipschitz condition to the notion of Lipschitz continuity.
\begin{itemize}
\item Show that any $L$-Lipschitz continuous vector field also satisfies \eqref{eq:osl} for a $\nu\geq L$.
\item Find an example of a vector field which satisfies \eqref{eq:osl}, but which is not Lipschitz continuous.
\end{itemize}
\end{exercise}
\last{There are several ways to model non-expansive and contractive vector fields. This work focuses on negative gradient flows, which we now introduce.}

\subsection{\last{Negative gradient flows}}
\label{sec:gradflow}
We now focus on a particular class of dynamical systems for which, by results in convex analysis, it is relatively immediate to verify the properties we have just derived. Let us consider a convex continuously differentiable function $V:\R^d\to \R$. A possible way to characterise the convexity of $V$ is through its gradient via the inequality
\[
\langle \nabla V(\y)-\nabla V(\x),\y-\x\rangle \geq 0,
\]
which is valid for every pair $\x,\y\in\R^d$. This condition implies that $X(\x)=-\nabla V(\x)$ is a non-expansive vector field. One could equivalently verify this property using \eqref{eq:lmax}, since the Hessian of a convex function is symmetric positive semi-definite, and hence $\lambda_{\max}(\partial_{\x}X(\x))=\lambda_{\max}(-\partial_{\x\x}^2 V(\x))\leq 0$ for all $\x,\y\in\R^d$, and hence $X(\x)=-\nabla V(\x)$ would be contractive by the reasoning from the previous section. 

The concept of $L$-smoothness from convex analysis is of particular importance to us for the applications in Section \ref{sec:examples}, since this is what will allow us to derive step size constraints for numerical discretisations of non-expansive flows. A convex, continuously differentiable, $V:\R^d\to\R$ is said to be $L$-smooth if its gradient is $L$-Lipschitz: 
\begin{equation}\|\nabla V(\x) - \nabla V(\y) \|_2\leq L \|\x - \y\|_2,
\label{eq:smooth_convex}
\end{equation}
for all $\x, \y\in\R^d$. \new{$L$-smoothness, convexity, and continuous differentiability are common assumptions in studies dealing with the convergence properties of gradient descent schemes, i.e., iterative schemes of the form $\x_{k+1}=\x_k - h_k\nabla V(\x_k)$.} An important result in convex analysis, the so-called \emph{Baillon--Haddad theorem}, tells us that this holds if and only if the following inequality holds for all $\x,\y\in\R^d$:
\begin{equation}
\langle \nabla V(\x) - \nabla V(\y), \x -\y\rangle \geq \frac{1}{L} \|\nabla V(\x) - \nabla V(\y)\|_2^2.
\label{eq:cocoercivity}
\end{equation}
\begin{exercise}
Assume that $V:\R^d\to\R$ is continuously differentiable, convex and $L$-smooth for some $L>0$. Prove the inequality in \eqref{eq:cocoercivity} using the following steps:
\begin{itemize}
    \item Use the fundamental theorem of calculus, applied to the scalar function $\varphi:[0,1]\to \R$ with $\varphi(t) = \nabla V(t\x + (1 - t) \mathbf z)$ to show that \eqref{eq:smooth_convex} implies the following inequality:
    \begin{equation}
    V(\mathbf z)\leq V(\x) + \langle\nabla V(\x), \x - \*z\rangle + \frac{L}{2}\|\x - \*z\|_2^2. \label{eq:smooth_convex_alternative}
    \end{equation}
    \item Add $\langle \nabla V(\y), \*z - \x\rangle$ to both sides of \eqref{eq:smooth_convex_alternative}, and minimise the left hand side with respect to $\mathbf z$ to obtain
    \[V(\y) + \langle\nabla V(\y), \x -\y\rangle \leq V(\x) +\langle \nabla V(\x) - V(\y), \x - \*z\rangle + \frac{L}{2} \|\*z - \x\|_2^2.\]
    \item Minimise the right hand side with respect to $\*z$, and add the resulting inequality to the corresponding inequality with $\x$ and $\y$ swapped. Upon rearranging, you should find \eqref{eq:cocoercivity}.
\end{itemize}
\end{exercise}
\last{Following the procedure in Section \ref{se:numericalNetworks}, we can build networks with layers based on negative gradient flows. Still, if we want those layers to be $1$-Lipschitz, we need to be careful when discretising, as described in the following section.}
\subsection{\last{$1$-Lipschitz networks based on gradient flows}}
Let us first consider a simple example to comment on the numerical approximation of the solutions of these non-expansive dynamical systems. Let $V(\x)=\|\new{\x}\|_2^2/2$. This is an $L$-smooth potential with $L=1$. The vector field $X(\x)=-\nabla V(\x)=-\x$ is hence non-expansive. We now suppose not to be able to exactly solve the system of differential equations $\dx(t)=X(\x(t))=-\x(t)$ and try to approximate its solutions at time $h>0$ with the explicit Euler method. This procedure leads to
\[
\x\mapsto \psi_X^h(\x):=\x - h \x = (1-h)\x,
\]
which provides an approximation of $\phi^h_{X}(\x)=e^{-h}\x$. We see that
\[
\left\|\psi_X^h(\y)-\psi_X^h(\x)\right\|_2 = |1-h|\cdot\|\y-\x\|_2\new{,\,\,\forall \x,\y\in\R^d},
\]
which is smaller than or equal to the initial distance $\|\y-\x\|_2$ if and only if $0\leq h \leq 2$. Therefore, even though we are considering one of the simplest dynamical systems, we see that we can not allow for an arbitrarily large time step $h$ if we want to numerically reproduce the non-expansiveness of the continuous solution $\phi^h_X$.

This reasoning extends to $L$-smooth convex potentials $V:\R^d\to\R$. Indeed, if we take steps using the Euler method with some step size $h>0$, $\new{\x}\mapsto \psi_{-\nabla V}^h(\x) = \new{\x} - h\nabla V(\x)$, we find that\new{, for every $\x,\y\in\R^d$,}
\begin{align*}
\|\psi_{-\nabla V}^h (\x) - \psi_{-\nabla V}^h(\y)\|_2^2 &= \|\x - \y -h(\nabla V(\x) - \nabla V(\y))\|_2^2\\
&= \|\x - \y\|_2^2 -2h\langle \nabla V(\x) - \nabla V(\y), \x - \y\rangle \\
&\qquad\qquad\qquad\qquad\qquad\qquad+h^2 \|V(\x) - V(\y)\|_2^2\\
&\leq \|\x - \y\|_2^2 + \Big(h^2- \frac{2h}{L}\Big) \|\nabla V(\x) - \nabla V(\y)\|_2^2,
\end{align*}
where the inequality follows directly from \eqref{eq:cocoercivity}. As a result, the Euler method preserves the non-expansiveness of the flow, as long as the step size $h$ satisfies the constraint $0\leq h \leq \frac{2}{L}$.

An example of a potential satisfying the requirements above is $V(\x)=\boldsymbol{1}^{\top}\gamma(A\x+\mathbf{b})$, where 
\[\gamma(\x)_i=\begin{cases}\frac{1}{2}x_i^2,\quad&\text{if}\quad x_i > 0,\\
0,\quad&\text{otherwise},\end{cases}\]
$\boldsymbol{1}\in\R^H$ is a vector of ones, and $A\in\R^{H\times d},\mathbf{b}\in\R^H$ are trainable weights. An Euler step with $X(\x)=-\nabla V(\x)$ leads to the layer
\begin{equation}\label{eq:gradientLayer}
\x\mapsto \psi^h_X(\x):=\x -h A^{\top}\sigma(A\x+\mathbf{b}),
\end{equation}
where $\sigma(\x)_i=\max\{0,x_i\}$ is the ReLU activation function. 
\begin{exercise}\label{ex:boundLipschitz}
Show that the vector field $X=-\nabla V$ in \eqref{eq:gradientLayer} is $L$-Lipschitz with $L=\|A\|_2^2$.
\end{exercise}
Based on the previous exercise, we can conclude that the Euler step in \eqref{eq:gradientLayer} is $1$-Lipschitz if 
\begin{equation}
0\leq h\leq 2/\|A\|_2^2.
\label{eq:stepsize_constraint}
\end{equation}
As we will see in Section \ref{sec:examples}, we can easily satisfy this constraint during training, using the power method to keep track of $\|A\|_2$. Of course, since non-expansive maps compose to give non-expansive maps, we can now stack any number of such blocks to get a non-expansive map.

%% file: hamiltonianNetworks.tex
\section{{Hamiltonian neural networks}}
\label{sec:hamiltonian}

In this section, we are revisiting the problem of stable training as outlined in Section \ref{sec:intro}. We now present a strategy to do so based on designing the network layers so they approximate the solution of some suitably constructed dynamical system. This idea is introduced and developed in \cite{galimberti2023hamiltonian}.

The dynamical systems we consider are canonical Hamiltonian systems. We define them on $\mathbb{R}^{2d}$ via the differential equations
\begin{equation}\label{eq:hamiltonian}
    \begin{cases}
        \dx(t)=\J\nabla H\left(\x(t)\right)=:X_H\left(\x(t)\right), \\
        \x(0)=\x_0,
    \end{cases}\,\,\J=\begin{pmatrix}
        0 & I \\ -I & 0
    \end{pmatrix}\in\R^{2d\times 2d},
\end{equation}
where $0,I\in\R^{d\times d}$ are the zero and identity matrices respectively, and $H\in\mathcal{C}^2\left(\R^{2d},\R\right)$ is called the Hamiltonian of the system. The so-called canonical symplectic matrix $\J$ is skew-symmetric, i.e., $\J^{\top}=-\J$. This structure of $\J$ implies that the energy function $H$ is conserved along the solutions of \eqref{eq:hamiltonian}:
\[
    \frac{\dd}{\dd t}H\left(\x(t)\right) = \langle\nabla H\left(\x(t)\right), \dx(t)\rangle =  \nabla H\left(\x(t)\right)^{\top}\J\nabla H\left(\x(t)\right) = 0.
\]
Another interesting property of Hamiltonian systems is that $\x(t)=\phi^t_{X_H}(\x_0)$ is a symplectic map for any $t\geq 0$, i.e.,
\begin{equation}\label{eq:symplectic}
    \partial_{\x_0} \x(t)^{\top}\J\partial_{\x_0} \x(t)=\J.
\end{equation}
\begin{exercise}[Proof of \eqref{eq:symplectic}]
The proof can be divided into the following two steps:
\begin{enumerate}
    \item Verify that \eqref{eq:symplectic} holds for $t=0$.
    \item Show that
    \[
    \frac{\dd}{\dd t}\left(\partial_{\x_0} \x(t)^{\top}\J\partial_{\x_0} \x(t)\right) = 0
    \]
    for any $t$. (\textbf{Hint:} Differentiate the system of Hamiltonian equations in \eqref{eq:hamiltonian} with respect to $\x_0$.)
\end{enumerate}
Proving these two points allows us to conclude since 
\[
\partial_{\x_0} \x(t)^{\top}\J\partial_{\x_0} \x(t)=\partial_{\x_0} \x(0)^{\top}\J\partial_{\x_0} \x(0) = J,
\]
as desired.
\end{exercise}
Equation \eqref{eq:symplectic} implies that
\[
    \left\|\J\right\|_2 = \left\|\partial_{\x_0} \x(t)^{\top}\J\partial_{\x_0} \x(t)\right\|_2 \leq \left\|\partial_{\x_0} \x(t)\right\|_2^2\|\J\|_2,
\]
and hence
\begin{equation}\label{eq:lowerBound}
    \left\|\partial_{\x_0} \x(t)\right\|_2\geq 1.
\end{equation}
There is a class of numerical methods, called symplectic, which allows to reproduce the property in \eqref{eq:symplectic} also on the approximate solution\new{, see \cite[Chapter VI]{hairer2011geometric} and \cite{leimkuhler2004simulating,Sanz-Serna_1992}}. A one-step method $\psi^h$ is symplectic if, when applied to a Hamiltonian system, it satisfies
\begin{equation}\label{eq:symplMethod}
    \left(\partial_{\x_0} \psi^h\left(\x_0\right)\right)^{\top}\J(\partial_{\x_0} \psi^h\left(\x_0\right))=\J.
\end{equation}
\begin{exercise}
Show that the composition $F\circ G:\R^{2d}\to\new{\R^{2d}}$ of two continuously differentiable symplectic maps $F,G:\mathbb{R}^{2d}\to\mathbb{R}^{2d}$ is again symplectic. We recall that a continuously differentiable map $F$ is symplectic if $\partial_{\x}F(\x)^\top J \partial_{\x}F(\x)=J$ for every $\x\in\R^{2d}$.
\end{exercise}

A Hamiltonian Neural Network (HNN) $\Phi$ is a network with $j$-th layer defined via a single step of a symplectic method $\psi^h$ applied to a parametrised Hamiltonian system with Hamiltonian function $H_j$. This construction removes the vanishing gradient problem since $\Phi$, being symplectic, satisfies \eqref{eq:normIneq}. We now conclude this section by providing an explicit example of an HNN.

Theoretically, there is no constraint on how the parametric Hamiltonian functions should be defined. However, some choices might restrict the expressiveness of the network or lead to network architectures completely different from the ones people are used to. A choice for $H_j$ that allows to recover expressive and commonly used architectures is
\begin{equation}\label{eq:parHamiltonian}
    H_j\left(\x\right)=\langle \boldsymbol{1},\gamma\left(A_j\x+\*a_j\right)\rangle,\,\,A_j\in\R^{2d\times 2d},\,\*a_j\in\R^{2d},
\end{equation}
where $\gamma:\R\to\R$ is a differentiable function applied to the entries of its input vector, and $\boldsymbol{1}\in\R^{2d}$ is a vector of all ones. This parametrisation allows us to get
\begin{equation}\label{eq:rhs}
    \J\nabla H_j\left(\x\right) = \J A^{\top}_j\sigma\left(A_j\x+\*a_j\right),
\end{equation}
where $\gamma'=\sigma$ becomes the activation function of the neural network. After having defined this parametric vector-valued function, one simple option is to define the network layers $\Phi_j$ as explicit Euler steps applied to vector fields as in \eqref{eq:rhs} to get
\begin{equation}\label{eq:euler}
    \Phi_{j}\left(\x\right) = \x + h\J A_j^{\top}\sigma\left(A_j\x+\*a_j\right),
\end{equation}
which has a similar structure to common ResNets. For example, to get $\sigma=\tanh$, one could set $\gamma=\log\circ \cosh$. The potential issue with defining the layer maps as in \eqref{eq:euler} is that the explicit Euler method is not symplectic and hence \eqref{eq:symplMethod} is not guaranteed. These considerations imply that even though we started from Hamiltonian systems for which \eqref{eq:lowerBound} holds, it might not be true that
\[
    \Big\|\partial_{\x} \Phi_j\left(\x\right)\Big\|_2\geq 1,
\]
and hence, there might still be a vanishing gradient problem. A solution to this issue is provided, for example, by the symplectic Euler method. Let us consider a splitting of the variable $\x\in\R^{2d}$ as $\x=(\q,\p)$, $\q,\p\in\R^{d}$. If the Hamiltonian function $H$ is separable, meaning that $H:\R^d\times\R^d\to\R$ is defined based on two functions $K,U:\R^d\to\R$ as $H\left(\q,\p\right)=K\left(\p\right)+U\left(\q\right)$, then the Hamiltonian dynamical system associated to $H$ writes
\[
    \begin{cases}
        \dot{\q}(t) = \nabla K\left(\p(t)\right),  \\
        \dot{\p}(t) = -\nabla U\left(\q(t)\right), \\
        \q(0)=\q_0,\,\p(0)=\p_0.
    \end{cases}
\]
The symplectic Euler method for this problem is explicit and reads
\begin{equation}\label{eq:symplEuler}
    \psi^h(\q,\p) =
    \begin{pmatrix}
        \widehat{\q} \\
        \p-h\nabla U\left(\widehat{\q}\right)
    \end{pmatrix},\quad \widehat{\q}=\q+h\nabla K(\p).
\end{equation}

\begin{exercise}
    Prove that the map $\psi^h$ in \eqref{eq:symplEuler} is symplectic, i.e., satisfies \eqref{eq:symplMethod}. (\textbf{Hint:} Write it as the composition of two simpler symplectic maps looking at how $\widehat{\q}$ is defined.)
\end{exercise}

In order to make the parametric Hamiltonian in \eqref{eq:parHamiltonian} separable, we can assume $A_j \in\R^{2d\times 2d}$ has a block structure as
\[
    A_j = \begin{pmatrix} 0 & B_j \\ C_j & 0 \end{pmatrix},\,\,B_j,C_j\in\R^{d\times d},
\]
and we also write $\*a_j=\begin{pmatrix} \*b_j^{\top} &  \*c_j^{\top}\end{pmatrix}^{\top}$, $\*b_j,\*c_j\in\R^d$. In this way, using the same partitioning $\x=\left(\q,\p\right)$ as before, we get
\[
    H_j\left(\q,\p\right)  =  
    \langle \boldsymbol{1},\gamma\left(B_j\p + \*b_j\right)\rangle +\langle\boldsymbol{1}, \gamma\left(C_j\q + \*c_j\right)\rangle=:K_j\left(\p\right) + U_j\left(\q\right),
\]
where $\boldsymbol{1}\in\R^d$ is the vector with all components equal to 1.
To conclude, we can then get an explicitly defined HNN with $j$-th layer
\begin{equation}\label{eq:HNNlayer}
    \psi_j(\x) =
    \begin{pmatrix}
        \widehat{\q} \\
        \p-hC_j^{\top}\sigma\left(C_j\widehat{\q}+\*c_j\right)
    \end{pmatrix},\quad \widehat{\q}:=\q+h B_j^{\top}\sigma\left(B_j\p+\*b_j\right),
\end{equation}
which does not suffer from vanishing gradient problems.

\begin{figure}[!htb]
\centering
\includegraphics[scale=1]{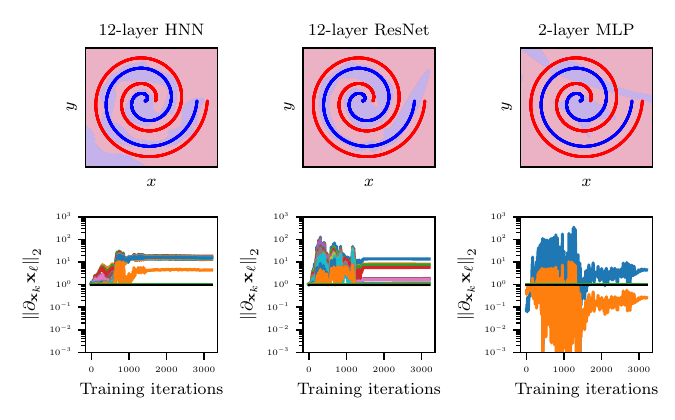}
\caption{A comparison of a 12-layer HNN, a 12-layer ResNet and a 2-layer MLP on the ``Swiss roll'' dataset, as previously considered for a 12-layer MLP in Figure \ref{fig:ex:training}. Both the HNN and ResNet attain a test accuracy of 100\%, while the 2-layer MLP has a test accuracy of 79.23\%. Note that the Jacobian norms behave much less extremely than they did for the 12-layer MLP in Figure \ref{fig:ex:training}, resulting in networks that train better.}
\label{fig:HNNexperiment}
\end{figure}

In the remaining part of this section, we provide a numerical experiment testing out the architectures we have derived and showing the improvements in terms of vanishing gradient issues provided by HNNs. We consider the problem of classifying into two classes the points in the 2D ``Swiss roll'' dataset, which can be seen in the top row of Figure \ref{fig:HNNexperiment}. The red and blue colours in the figure represent the two classes. We test different network architectures. The first is the HNN with layers defined as in \eqref{eq:HNNlayer}, the second is a ResNet with layers based on the explicit Euler method and of the form $\Phi_j(\x)=\x + hB_j^{\top}\sigma(A_j\x+\*{b}_j)$, and the third is an MLP with layers defined as $\Phi_j(\x)=B_j^{\top}\sigma(A_j\x + \*{b}_j)$. We consider the HNN and ResNet with $\ell=12$ hidden layers of the form above (as we did for the MLP in Figure \ref{fig:ex:training}), composed with a final linear layer to adapt the network to the output dimensionality which, in this case, is two. The MLP is also considered for the case of $\ell=2$ hidden layers. The dataset is embedded in a higher-dimensional space of dimension four in the following way $(x_1, x_2) \mapsto (x_1, 0, x_2, 0)$. The network layers then preserve this intermediate fixed dimension.

In Figure \ref{fig:HNNexperiment}, we can see that the ResNet and HNN models both perform accurately on this simple task, leading to a $100\%$ classification accuracy over a test set. Instead, the MLP with $12$ layers does not train appropriately, as we saw in Figure \ref{fig:ex:training}, leading to a classification that is only slightly better than chance. On the other hand, the MLP with two layers trains slightly better, leading to around $80\%$ accuracy. These four models have been chosen to illustrate the issue of having vanishing gradients and, consequently, not being able to train the network. In the bottom row of Figure \ref{fig:HNNexperiment}, we plot the norms of the Jacobian matrices of the last hidden layer with respect to the previous ones throughout the training iterations. For each of the four models, a fixed test data point has been used to evaluate these Jacobian matrices. We see that the ResNet and HNN models lead to well-behaved Jacobians. On the other hand, the MLP model has vanishing gradient issues, which lead to the impossibility of training the model with $12$ layers, whereas these issues do not arise when training a network with just two hidden layers. While the HNN is built so that the norm of the Jacobian is never smaller than one, as can be seen in the plot, the skip connections in the ResNet naturally lead to stable behaviour under suitable weight initialisation. This is not surprising since residual connections were introduced precisely to allow the training of deeper networks.

%% file: asymptoticStability.tex
\section{Networks with stable equilibria}
\label{sec:stableEquilibria}

Up to now, we have seen how the stability of dynamical systems can be characterised either in terms of the reciprocal behaviour of pairs of trajectories, in Section \ref{sec:nonexpansiveness}, or in terms of the presence of conserved quantities, as in Section \ref{sec:hamiltonian}. Another typical way to analyse the stability of dynamical systems is through their stationary points, also called equilibria. Let us consider the time-independent dynamical system described by the system of differential equations $\dx(t) = X(\x(t))$, for some right-hand side $X:\mathbb{R}^d\to\mathbb{R}^d$. The equilibria of this system define the set $\mathcal{E}=\left\{\bar{\x}\in\R^d:\,X(\bar{\x})=0\right\}$. The peculiarity of these points is that the solutions of the differential equation with initial conditions in $\mathcal{E}$ will be trivial. More explicitly, $\phi^t_X(\x_0)=\x_0$ for every $t\geq 0$ when $\x_0\in\mathcal{E}$. The points in $\mathcal{E}$ can be characterised in terms of their stability properties, as we formalise in the following definition.
\begin{definition}
Let $\bar{\x}\in\mathcal{E}$ be an equilibrium point of the system of differential equations $\dx(t)=X(\x(t))$. We say $\bar{\x}$
\begin{itemize}
    \item \textit{stable} if for every $\varepsilon>0$, there exists a $\delta>0$ such that if $\left\|\x_0-\bar{\x}\right\|_2<\delta$, it follows $\left\|\phi^t_X(\x_0)-\bar{\x}\right\|_2<\varepsilon$ for all $t\geq 0$,
    \item \textit{locally asymptotically stable} if there is a neighbourhood $B_{\new{\bar{\x}}}\subset\R^d$ of $\new{\bar{\x}}$ such that $\lim_{t\to+\infty}\left\|\phi^t_X(\x_0)-\bar{\x}\right\|=0$ whenever $\x_0\in B_{\bar{\x}}$,
    \item \textit{globally asymptotically stable} if $\lim_{t\to+\infty}\left\|\phi^t_X(\x_0)-\bar{\x}\right\|_2= 0$ for every $\x_0\in\R^d$.
\end{itemize}
\end{definition}
\begin{figure}[!htb]
\centering
\includegraphics[scale=1]{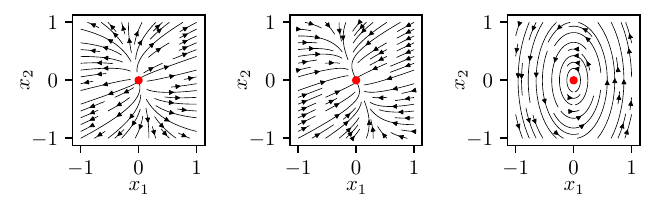}
\caption{Phase portrait of three linear systems. From left to right, the origin is asymptotically stable, unstable, and stable but not asymptotically.}
\label{fig:stability}
\end{figure}
\begin{exercise}
Suppose that the matrix $A\in\R^{d\times d}$ is diagonalisable. 
\begin{itemize}
    \item Prove that if the real parts of the eigenvalues of $A$ are all strictly negative, then the dynamical system $\dx(t)=A\x(t)$ has a unique equilibrium point at the origin, and it is globally asymptotically stable. 
    \item Consider again $\dx(t)=A\x(t)$. What is a condition on the eigenvalues of $A$ leading to a system which is stable but not asymptotically stable? Find some examples of matrices satisfying this condition.
\end{itemize}
\end{exercise}
We plot in Figure \ref{fig:stability} the phase portrait of three linear dynamical systems for which the origin has different stability properties. Thicker lines correspond to faster dynamics, meaning that the norm of $X$ is bigger.

The study of the stability of equilibria is a very well-developed area in the field of dynamical systems. One of the most commonly used tools to study their stability is the notion of Lyapunov function. 
\begin{definition}[Lyapunov function]
A continuously differentiable function $V:U\to \R$, $U\subset \R^d$ open, is a Lyapunov function for $\dx(t)=X(\x(t))$ associated to the equilibrium point $\bar{\x}\in U$ if it satisfies
\begin{itemize}
    \item $V(\x)> 0$ for every $\x\in U\setminus\{\bar{\x}\}$, and $V(\bar{\x})=0$,
    \item $\frac{\dd}{\dd t}V(\x(t))=\langle\nabla V(\x(t)), \dx(t)\rangle = \langle\nabla V(\x(t)), X(\x(t))\rangle\leq 0$ for every solution curve $t\mapsto \x(t)$ with $\x(0)\in U$.
\end{itemize}
We call $V$ a strict Lyapunov function if $\frac{\dd}{\dd t}V(\x(t))<0$.
\end{definition}
Geometrically, Lyapunov functions lead to subsets of $\mathbb{R}^d$ from which the solution can not escape. These subsets are the sub-level sets of $V$, defined as $L_c = \left\{\x\in\mathbb{R}^d:\,\,V(\x)\leq c\right\}$. Indeed, since the gradient vector field $\nabla V(\x)$ is orthogonal to the level sets of $V$, the condition $\langle\nabla V(\x), X(\x)\rangle\leq 0$ corresponds to saying that the vector field $X$ does not point outward of the sublevel sets. In other words, the value taken by a Lyapunov function at time $0$, $V(\x_0)=c$, has to be an upper bound of the value at any time $t\geq 0$, i.e., $V(\phi^t_X(\x_0))\leq c$ meaning that $\phi^t_X(\x_0)\in L_c$ for every $t\geq 0$. Based on this reasoning, we can conclude that the presence of a Lyapunov function guarantees the stability of the associated equilibrium point $\new{\bar{\x}}$. 
\begin{exercise}
Find a Lyapunov function for the system
\[
\begin{cases}
    x' = -x + y^2,\\
    y' = -2y+3x^2.
\end{cases}
\]
(\textbf{Hint:} Consider $V(x,y)=ax^2+bxy+cy^2$ for a suitable choice of $a,b,c\in\R$.)
\end{exercise}
\begin{exercise}
Show that if $\dx(t)=X(\x(t))$ admits a strict Lyapunov function $V:U\to\R$, for an equilibrium point $\bar{\x}\in U$, then $\bar{\x}$ is locally asymptotically stable.
\end{exercise}

\subsection{Learning stable dynamical systems}
To describe the time evolution of physical systems, one needs governing equations, specifically the right-hand side $X:\R^d\to\R^d$ of a differential equation. Traditionally, experts in the field have created these models by deriving an accurate description of the system. However, with modern computational power and abundant data, data-driven modelling is gaining considerable attention. When the system's behaviour is partially known (e.g., it has a stable equilibrium point), the approximate model should reflect these properties.

In \cite{kolter2019learning}, the authors explicitly build a data-driven model $X:\R^d\to\R^d$ which is known to have a Lyapunov function $V:\R^d\to\R$ associated to a prescribed equilibrium point $\bar{\x}\in\R^d$. To do so, they parametrise $X$ as follows
\begin{equation}\label{eq:stableVecField}
X(\x) = \hat{X}(\x) - \nabla V(\x)\frac{\mathrm{ReLU}\left(\nabla V(\x)^{\top}\hat{X}(\x)+\mu V(\x)\right)}{\left\|\nabla V(\new{\x})\right\|_2^2},\,\,\mu>0,
\end{equation}
where $\hat{X}:\R^d\to\R^d$ can be any neural network, while $V:\R^d\to \R$ is modelled as a positive-definite scalar-valued neural network which is guaranteed to be convex in the inputs and has the correct stationary point. 
\begin{exercise}
Prove that $V:\R^d\to\R$ is a Lyapunov function for $X$ in \eqref{eq:stableVecField}.
\end{exercise}
In \cite{jaffeLearningNeuralContracting2024}, the authors propose the parametrisation $X(\x)=A(\x,\bar{\x})(\x-\bar{\x})$ with 
\begin{equation}\label{eq:symND}
\mathrm{Sym}(A(\x,\new{\bar{\x}}))=\frac{1}{2}\left(A(\x,\bar{\x})+A(\x,\bar{\x})^\top\right),
\end{equation}
which is negative definite. In their case, $A:\R^d\times\R^d\to\R^{d\times d}$ is a matrix valued neural network forced to satisfy \eqref{eq:symND}, and $X(\bar{\x})=0$ so $\bar{\x}\in\R^d$ is an equilibrium point of $X$. Theorem 3 in \cite{jaffeLearningNeuralContracting2024} shows that this parametrisation for $X$ ensures the asymptotic stability of $\bar{\x}$.

\subsection{Asymptotic stability for adversarial robustness}
We have already seen in Section \ref{sec:nonexpansiveness} how building neural networks based on non-expansive dynamical systems can improve their robustness to input perturbations and, hence, also to adversarial attacks. Contractive dynamical systems also have an asymptotically stable equilibrium point, so one might want to investigate how directly focusing on networks based on asymptotically stable dynamical systems can be beneficial in this context. This point of view has been considered in several works, like \cite{huang2022fi,kang2021stable}. In these papers, the authors exploit that with (locally) asymptotically stable equilibria, the trajectories starting in an open neighbourhood of the equilibria converge to the equilibrium to ensure that the network prediction is not too sensitive to input perturbations.

%% file: worked_examples.tex
\section{Worked examples}
\label{sec:examples}
Let us now dive in and work through two applications of stable neural networks!

We will develop two examples demonstrating the implementation and results of non-expansive neural networks applied to solving an ill-conditioned two-dimensional inverse problem and classifying images robustly. This section includes a description of the two examples, together with the results one can get running the two associated \texttt{jupyter} notebooks in the repository related to this paper\footnote{\href{https://github.com/davidemurari/bookChapterDS}{https://github.com/davidemurari/bookChapterDS}}. The notebooks are called \texttt{inverse\_problem.ipynb} and \texttt{adversarial\_robustness.ipynb}. This section and the notebooks are meant to be self-contained.

We implement a non-expansive neural network following the principles presented in Section \ref{sec:nonexpansiveness}. This will be used for both examples. Each network layer corresponds to an explicit Euler step of a suitable vector field, defined in \new{the code snippet} \ref{snippet:nonexpansive_block}. This block is based on linear layers. We extend it to convolutional layers in the notebooks associated with this paper. To create an object of the \texttt{NonExpansiveBlock} class, we need to specify the input dimension, the output dimension, and the final time of the numerical integration. The forward method takes as input the current position and the number of substeps we want to take to reach the final time and it returns the updated position.  

%
\begin{code}
\footnotesize
\captionof{listing}{Neural network block based on numerically integrating a non-expansive ODE\new{, with layer $\x\mapsto \x - hA^\top \mathrm{ReLU}(A\x+\mathbf{b})$, $h=T/\texttt{n\_steps}$.}}
\begin{lstlisting}[style=pyLikeMinted]
class NonExpansiveBlock(torch.nn.Module):
    def __init__(self, dim_inner=10, dim_outer=10, T=1.):
        super().__init__()
        self.lin = torch.nn.Linear(dim_inner, dim_outer)
        self.T = T
        
    def forward(self, x, n_steps):
        for i in range(n_steps):
            x = x - (self.T / n_steps) * relu(self.lin(x)) @ self.lin.weight.T
        return x
\end{lstlisting}
\label{snippet:nonexpansive_block}
\end{code}

A non-expansive neural network can be obtained by composing several of these blocks. For this purpose, we need to implement the step size constraint from \eqref{eq:stepsize_constraint}, which requires us to estimate the spectral norms of the linear layers. This is done with the power method. Let us consider a matrix $A\in\R^{d\times c}$ defining the linear layer of interest. The power method is implemented as
\new{\begin{equation}\label{eq:powerMethod}
\new{\bu}_{i+1} =\frac{A^\top A\new{\bu}_i}{\|A^\top A\new{\bu}_i\|_2},\,\,\,i=0,\ldots,k-1,
\end{equation}
with $\new{\bu_0}\in\R^c$.} The vector $\new{\bu_0}$ could either be an initial estimate of the first right singular vector of $A$ or a random vector. If $\new{\bu_0}$ is not orthogonal to the target right singular vector, this iteration computes $\new{\bu_k}$ which usually approximates the first right singular vector of $A$, and $\sqrt{\|A^\top A\new{\bu_k}\|_2}$ converges to $\|A\|_2$ as $k\to \infty$.
\begin{exercise}
Implement the power method as described in \eqref{eq:powerMethod} and verify that it provides an accurate approximation of the spectral norm of the following matrices:
\begin{itemize}
    \item $A=I$, i.e., the $10\times 10$ identity matrix,
    \item $A=5 I$,
    \item $A = \exp(B-B^\top)$ for a random matrix $B\in\R^{10\times 10}$, with $\exp$ being the matrix exponential in this context. This is an orthogonal matrix, i.e., $A^\top A=I$, so we would expect it to have norm $1$.
\end{itemize}
\end{exercise}
Before training, we run many iterations of the power method and save the resulting estimates of the top right singular vectors of the linear layers. Much of the usual training loop for neural networks remains the same for networks built using non-expansive blocks: we load a minibatch of data and pass it through the network, we evaluate the network's predictions using the loss function, and we backpropagate and perform a gradient update. Before passing to the next minibatch, however, we update our estimates of the spectral norms of the weights using the power method. Since we have a good estimate of the top right singular vector, we use it to warm-start the power method, making it possible to use just a single iteration of the power method. After this, \texttt{n\_steps} in \texttt{forward} in snippet \ref{snippet:nonexpansive_block} is computed as the smallest integer $n$ such that $h =T/n$ (with $T$ the total integration time) satisfies the step size constraint in \eqref{eq:stepsize_constraint}. That is to say, we adapt the step size as necessary to preserve non-expansiveness when the Lipschitz constant of the vector field grows.

In the image classification example, we also compare the non-expansive network to a baseline ResNet. The proposed implementation follows a structure similar to the non-expansive network. The main change is in \texttt{ResidualBlock}, which implements the explicit Euler step of a different differential equation, which does not generally have a non-expansive flow. We present the residual block in snippet \ref{snippet:resnet}. Again, this block is described in the case of linear layers for simplicity, but the implementation in the notebooks is extended to convolutional layers.

%
\begin{code}
\footnotesize
\captionof{listing}{ResNet architecture\new{, with layer $\x\mapsto \x + B\mathrm{ReLU}(A\x+\mathbf{b})$.}}
\begin{lstlisting}[style=pyLikeMinted]
class ResidualBlock(nn.Module):
    def __init__(self, dim_inner=10, dim_hidden=10):
        super().__init__()
        self.linearA = torch.nn.Linear(dim_inner,dim_hidden)
        self.linearB = torch.nn.Linear(dim_hidden,dim_inner)

    def forward(self, x):
        return x + self.linearB(relu(self.linearA(x)))
\end{lstlisting}
\label{snippet:resnet}
\end{code}

\subsection{Ill-conditioned inverse problem}
\label{sec:inverse_problem}
Recall the inverse problem shown in Section \ref{sec:intro}: we are tasked with inverting measurements $\*y^\delta\in \R^m$ taken of some ground truth vector $\*x^\dagger$, where
\begin{equation}
\*y^{\delta} = A\*x^\dagger + \*z.
\label{eq:inverse_problem}
\end{equation}
For concreteness, we will here consider the simple case where the ground truth vectors $\x$ are supported on a curved set in $\R^2$, as shown in Figure \ref{fig:ip}, and its forward operator is given by
\[
    A = \begin{pmatrix} 1+\varepsilon & 1\\
1 & 1+\varepsilon\end{pmatrix}.
\]
for $\varepsilon = 1/4$. Additionally, $\*z$ is given by Gaussian white noise, and we can use knowledge of $\*y^\delta$ and $A$ to estimate $\*x^\dagger$. In addition to the test set, which is shown in Figure~\ref{fig:ip}, we are given a training set of moderate size consisting of pairs of $\*x$ and matching, noisy, measurements $\*y^\delta$, which we may use to tune the parameters of any method under consideration. \new{That is to say, we are framing the problem of optimally regularising the inverse problem as a supervised learning problem. This can be contrasted with classical unsupervised approaches such as the Morozov discrepancy principle \cite{englRegularizationInverseProblems2000}. The Morozov discrepancy principle has the advantage of only using the measurements and an estimate of the noise level, but it is significantly less powerful than the supervised approach and is typically only used to tune a single parameter representing the regularisation strength.} The details of setting up the data are laid out in Section \enquote{Setting up the data} of the associated \texttt{jupyter} notebook, \texttt{inverse\_problem.ipynb}.
\begin{figure}[!htb]
\centering
\includegraphics[scale=1]{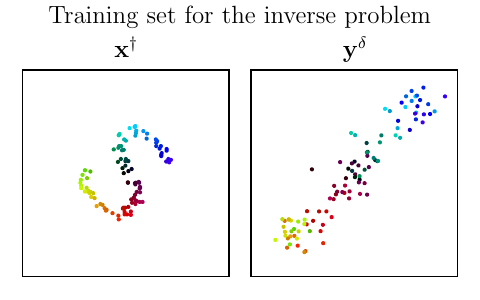}
\caption{The training set that we will use to tune methods for solving the inverse problem of determining $\*x^\dagger$ from $\*y^\delta$ in Equation \ref{eq:inverse_problem}.}
\label{fig:inverse_problem_training_set}
\end{figure}

\subsubsection{A classical regularisation approach}
\begin{figure}[!htb]
\centering
\includegraphics[scale=1]{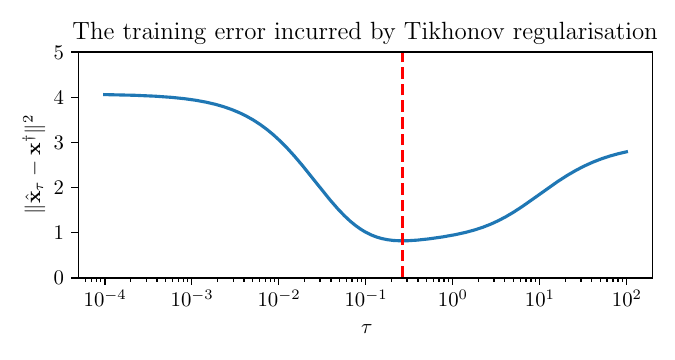}
\caption{We evaluate the performance of the reconstruction in Equation \ref{eq:tikhonov} (averaged over the training set) over a wide range of values of $\tau$ and select the value that attains the lowest error.}
\label{fig:tikhonov_tuning}
\end{figure}

Here, we will follow what is done in Section \enquote{Classical regularisation} of the associated \texttt{jupyter} notebook. As shown in Figure \ref{fig:ip}, na\"ively applying the inverse of $A$ is hopeless as a way to recover $\*x^\dagger$ from $\mathbf y^\delta$: the noise present in the measurement is blown up, obscuring any trace of the true signal. Classical approaches to overcoming this issue stabilise the inversion process by appropriately balancing the fit to measurements and fit to (some notion of) prior information. One of the most famous such regularisation methods is \emph{Tikhonov regularisation}, which introduces a regularisation parameter $\tau > 0$ to estimate $\*x$ by 
\begin{equation}
    \hat{\*x}_\tau = (A^\top A + \tau I)^{-1} A^\top \*y^\delta.
    \label{eq:tikhonov}
\end{equation}
Equivalently, $\hat{\*x}_\tau$ can be characterised as the unique minimiser of the functional $\*x\mapsto \|A\*x - \*y^\delta \|^2 + \tau \|\*x\|^2$, showing that this method naturally balances fitting the measurements (the first term), with ensuring that the estimate does not have a large norm (the second term). Since we have a training data set, and a 1-dimensional parameter $\tau$, we can optimise it by a simple (logarithmic) grid search, as shown in Figure \ref{fig:tikhonov_tuning}.

\begin{figure}[!htb]
    \centering
    \includegraphics[scale=1]{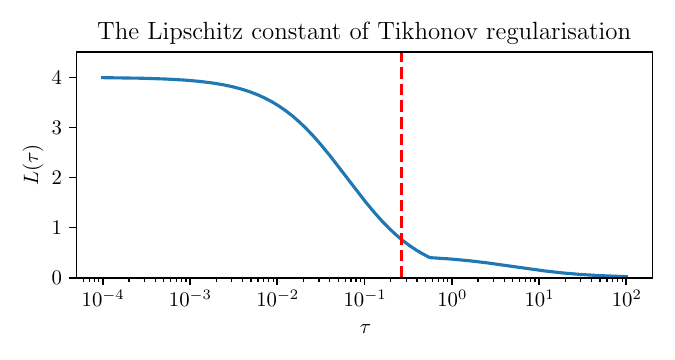}
    \caption{The Lipschitz constant of the reconstruction map for Tikhonov regularisation, seen as a function of the parameter $\tau$.}
    \label{fig:tikhonov_lipschitz}
\end{figure}

Although we have thought of this method as being parametrised by $\tau$, we can equivalently think of it as being parametrised by the Lipschitz constant of the corresponding reconstruction map
\[(A^\top A + \tau I)^{-1} A^\top,\]
which is just its operator norm, since it is a linear map. In fact, given the singular values $\{\sigma_i\}_i$ of $A$, this Lipschitz constant is given by
\begin{equation}
    L(\tau) = \max_i \frac{\sigma_i}{\sigma_i^2 + \tau},
    \label{eq:tikhonov_lipschitz}
\end{equation}
showing that it is monotonically decreasing in $\tau$, taking values between $0$ and $(1 / \min_i \sigma_i)$, as illustrated in Figure \ref{fig:tikhonov_lipschitz}.

\begin{figure}[!htb]
    \centering
    \includegraphics[scale=1]{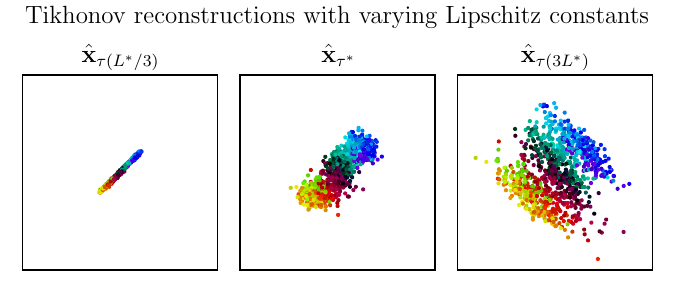}
    \caption{Test set reconstructions using the Tikhonov regularisation method, with the optimal parameter choice in the middle, more stable reconstructions on the left, and less stable reconstructions on the right.}
    \label{fig:tikhonov_recons}
\end{figure}

In particular, given the optimal parameter $\tau^*$, we can compute the corresponding Lipschitz constant $L^*:=L(\tau^*)$ and consider the behaviour of the reconstructions with Lipschitz constants $L^*/3$ and $3L^*$, say, corresponding to a more stable and less stable reconstruction than the optimal reconstruction, respectively. The results of doing this are shown in Figure \ref{fig:tikhonov_recons}. While the optimal parameter choice has stabilised the inversion, it is evident that neither it nor the other choices of the parameter $\tau$ allow for a faithful reconstruction of the curved shape of the support of the true data.
\subsubsection{Inversion using a stable neural network}
\begin{figure}[!htb]
    \centering
    \includegraphics[scale=1]{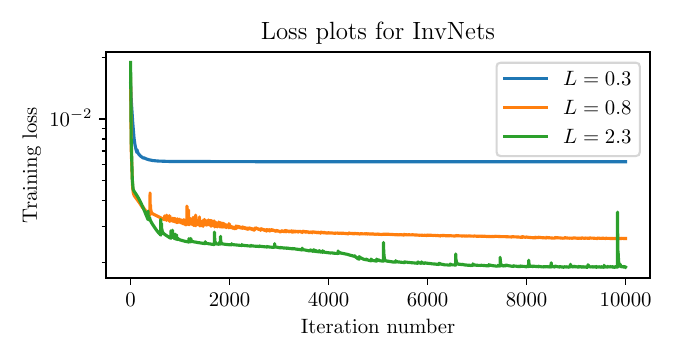}
    \caption{The evolution of the training loss function value over training for the three InvNets under consideration.}
    \label{fig:invnet_loss_plots}
\end{figure}

We will now overcome the shortcomings of Tikhonov regularisation, shown in Figure \ref{fig:tikhonov_recons}, by using a dynamical-systems-based neural network, which we will call InvNet. In addition to the learnable weights of this network, we will have a parameter $L > 0$, which serves as an upper bound on the Lipschitz constant of the network, and we will consider the choices $L=L^*/3$, $L=L^*$ and $L=3L^*$ as in Figure \ref{fig:tikhonov_recons}, with $L^*$ the Lipschitz constant of the optimal Tikhonov reconstruction map. An InvNet with choice of upper bound on Lipschitz constant $L$ will be denoted InvNet$_L$, and takes the following form, with each $\Phi_i$ a block of the form described in Snippet \ref{snippet:nonexpansive_block}:
\[\mathrm{InvNet}_L(\mathbf y) = c \cdot\mathrm{project}\circ \Phi_{\ell}\circ \cdots \circ \Phi_{1} \circ \mathrm{lift}(\mathbf y).\]
Here, $c$ is a learnable scalar parameter initialised to $1$, which is clamped between $-L$ and $L$, and the lifting and projection layers take a particularly simple: lift concatenates a vector of zeros to the input to fill out the dimensions, while project ignores the extra dimensions, so that both operations are clearly 1-Lipschitz. Finally, as described above, the composition of the dynamical blocks, $\Phi_\ell \circ \cdots \circ \Phi_1$, is kept non-expansive during training by keeping track of the operator norms of the weights and splitting the integration interval into more steps where necessary. As outlined in the Section ``A neural network approach'' of the \texttt{jupyter} notebook, we train three InvNets, InvNet$_{L^*/3}$, InvNet$_{L^*}$ and InvNet$_{3L^*}$. Concretely, we take $\ell=5$ dynamical blocks and lift from 2-dimensional inputs to 10-dimensional intermediate representations. For each network, we run 10,000 iterations of the Adam optimiser on the supervised loss function:
\[\frac{1}{n}\sum_{i=1}^{n}\|\mathrm{InvNet}_{L}(\y_i^\delta) - \x_i^\dagger\|^2.\]
\new{The optimiser uses the learning rate $10^{-3}$, no weight decay, and the default PyTorch parameters for this method.} As should be expected, the networks become less constrained with increasing $L$, corresponding to lower training losses, which is confirmed by the plots in Figure \ref{fig:invnet_loss_plots}.

With these choices, we find that the reconstructions of both InvNet$_{L^*}$ and InvNet$_{3L^*}$ are better than any of the reconstructions done using Tikhonov regularisation. In particular, the non-linear nature of the neural networks allows us to capture the curved shape of the underlying dataset. Interestingly, InvNet$_{3L^*}$ performs quite well, even though this Lipschitz constant corresponds to unstable reconstructions for Tikhonov regularisation (recall Figure \ref{fig:tikhonov_recons}). This may be explained by the fact that $3L^*$ is an upper bound for the Lipschitz constant of InvNet$_{3L^*}$, which need not be tight.
\begin{figure}[!htb]
    \centering
    \includegraphics[scale=1]{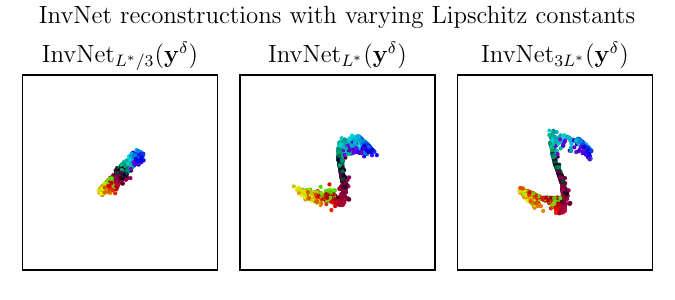}
    \caption{Test set reconstructions using InvNet, trained with three different upper bounds on the Lipschitz constant: $L=L^*/3$, $L=L^*$ and $L=3L^*$.}
    \label{fig:invnet_recons}
\end{figure}

\begin{exercise}
    Experiment with the notebook, exploring the following avenues:
    \begin{itemize}
    \item Vary the size of the training set: what happens with a much smaller training set of size $20$, or a much larger training set of size $500$? Comment on the effect of using a large $L$ as the size of the training set varies. \new{The notebook is set up so that one can input a single value for $L$ to facilitate this exercise.}
    \item The analysis in Section \ref{sec:gradflow} is not restricted to using ReLU as the activation function. Propose a different activation that works, meaning that it is Lipschitz continuous and non-decreasing, compute the associated step size constraint as in \eqref{eq:stepsize_constraint}, and implement the change in the block definition. How does the performance with the alternative activation function compare with using the $\mathrm{ReLU}$ activation function?
    \end{itemize}
\end{exercise}

\new{The example developed in this section is low-dimensional so that the simulations are faster, and the results are easier to visualise. However, the presented procedure can also be applied to higher-dimensional problems. Furthermore, the use of Lipschitz constraints in inverse problems has been very popular in the inverse problems literature. For example, the same architecture considered in this section was applied in \cite{sherryDesigningStableNeural2024} to design a provably convergent algorithm for image deblurring. Lipschitz networks in inverse problems can also be found in \cite{lunz2018adversarial,hasannasab2020parseval,hertrich2021convolutional,ryu2019plug}}.

\subsection{Adversarially robust image classification}
\label{sec:robust_image_classification}
\begin{figure}[!htb]
    \centering
    \includegraphics[scale=1]{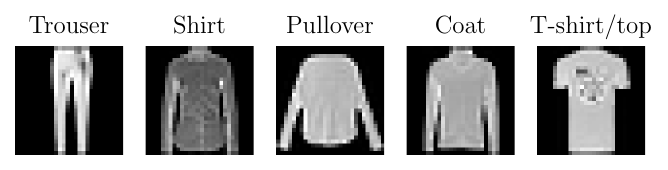}
    \caption{5 example images from the Fashion MNIST dataset.}
    \label{fig:sampleImages}
\end{figure}

In this section, we provide the details of the methods specific to this example and present the results one can obtain running the complete notebook \texttt{adversarial\_robustness.ipynb}. We work with the Fashion MNIST dataset, consisting of images of items from Zalando, along with a label denoting one of ten possible classes. It is based on a training set of 60,000 images and a test set of 10,000 images. Each is a greyscale image of size $28\times 28$. Figure~\ref{fig:sampleImages} shows five images in the training set with their associated labels. 
We implement a training routine as described above to train the neural network to classify the training images accurately, using the Adam optimiser and a one-cycle learning rate schedule. \new{The loss function is regularised using $\ell^2$-weight decay with penalty parameter $\gamma=10^{-5}$. We employ a one-cycle learning rate scheduler, starting from a minimum learning rate of $10^{-4}$ and peaking at $10^{-2}$, then annealing back to $10^{-4}$ via a cosine strategy over the total training steps.}

Once the network is trained, we can test its robustness to adversarial attacks. We consider the $\ell^2$-PGD attack, standing for Projected Gradient Descent based on the $\ell^2$-norm. The algorithm defining this attack is implemented in the notebook, but let us describe the mechanics of the attack in some more detail here. This attack aims to maximise the loss function \texttt{loss\_fn}, which we provide as input, by perturbing the input image \texttt{image}. The correct label for the input image is \texttt{target}, and the perturbation of the input image we allow has $\ell^2$-norm smaller than \texttt{epsilon}. To build this perturbation, we perform \texttt{n\_iter} iterations of the following procedure. Let us consider the function
\[
F(\texttt{delta}):=\texttt{loss\_fn}(\texttt{image} + \texttt{delta},\texttt{target}).
\]
Each of the \texttt{n\_iter} iterations consists of one step of size \texttt{step\_size} in the direction of $\partial_{\texttt{delta}}F(\texttt{delta})/\left\|\partial_{\texttt{delta}}F(\texttt{delta})\right\|_2$ followed by a projection over the $\ell^2$-ball of radius \texttt{epsilon} centered at the origin. Finally, \texttt{delta} is added to \texttt{image} to get the perturbed image.

\begin{figure}[!htb]
    \centering
    \includegraphics[scale=1]{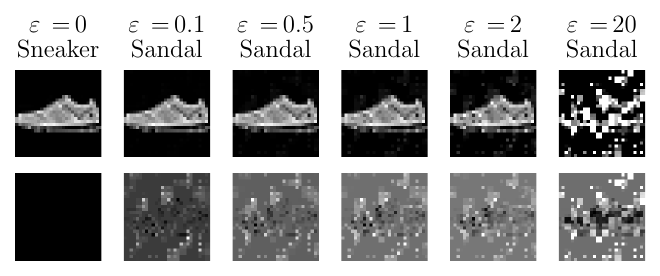}
    \caption{The first row displays the attacked images with increasing perturbation magnitudes. The second row displays the difference between the attacked and clean images. The titles specify the norm of the perturbation $\varepsilon$ and the ResNet prediction when given that image as an input.}
    \label{fig:attacks}
\end{figure}

In Figure \ref{fig:attacks} we show an example of an image attacked with the $\ell^2$-PGD attack with 100 iterations. The attack is displayed in different magnitudes, and one can see that the image looks increasingly distinct from the first one on the left, i.e., the clean image. The network we attack to obtain these perturbations is a ResNet trained to classify the test images with around $89\%$ accuracy.

\begin{figure}[!htb]
    \centering
    \includegraphics[scale=1]{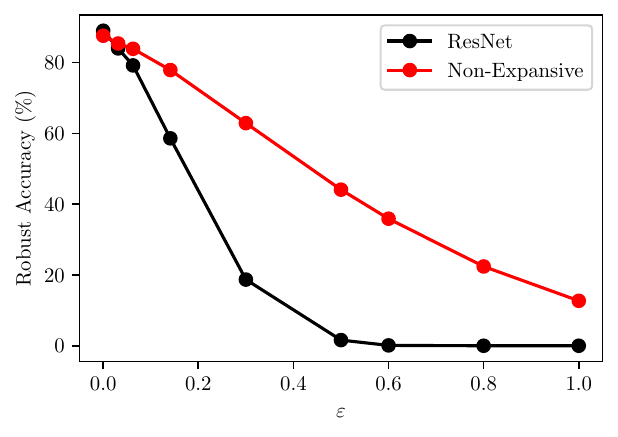}
    \caption{Comparison of the classification accuracy of a non-expansive network and a ResNet trained on 30,000 training images of the Fashion MNIST dataset. We then attack 1,000 of the test images with $100$ iterations of $\ell^2$-PGD of varying intensity $\varepsilon$. The attack magnitude is represented on the horizontal axis, with $\varepsilon=0$ corresponding to the clean images. The vertical axis displays the classification accuracy obtained with the attacked images. }
    \label{fig:robustAccuracy}
\end{figure}

We have now discussed all the necessary methods to evaluate the robustness of a non-expansive network architecture and compare it to that of an unconstrained ResNet. This comparison relies on two steps: training both networks on clean images and testing their accuracy on adversarial images built for the specific weights obtained after training. To have a code that takes five to ten minutes to train locally, we restrict the training and test sets to 30,000 and 1,000 images, respectively. The non-expansive network and the ResNet reach a similar test accuracy of around $88 \%$ to $89\%$. We train both models for $30$ epochs, again to benefit in terms of speed. When the training is completed, we freeze their parameters and build adversarial examples. The examples are obtained with 100 iterations of the $\ell^2$-PGD attack, and we generate them for different perturbation magnitudes. We consider eight attack magnitudes smaller than one and compare them with the clean accuracy corresponding to $\varepsilon=0$. Generating the attack for the 1,000 images takes around five to ten minutes. We plot the results obtained following this procedure in Figure \ref{fig:robustAccuracy}. We see a very small drop in performance for this relatively simple dataset when constraining the Euler steps to be $1$-Lipschitz. At the same time, robust accuracy improves over that of unconstrained layers. The gain in robustness is also expected for other datasets, while typically, the clean accuracy tends to decrease a bit more compared with the unconstrained model.

\begin{exercise}[Playing with the code]
Start from the \texttt{jupyter} notebook associated with this section, and test how the robustness changes by varying the \texttt{margin} in the loss function, the number of steps in the $\ell^2$-PGD attack, and the number of training epochs. Explore replacing the Fashion MNIST dataset and using other benchmark datasets, such as MNIST or CIFAR-10. If the training time grows considerably, we advise looking into \url{https://www.kaggle.com}, which offers $30$ free GPU hours per week. The code is already implemented to be accelerated with CUDA in case it is available on the machine.
\end{exercise}

\new{In this section, we focused on the Fashion MNIST dataset, which is fairly low-dimensional. The same design strategy for the network we considered was used to design adversarially robust networks for larger datasets such as CIFAR-10 or CIFAR-100 in \cite{celledoniDynamicalSystemsBased2023,sherryDesigningStableNeural2024,meunier2022dynamical}. Furthermore, Lipschitz networks different from those discussed in this work have been shown to be robust to adversarial examples (see, e.g., \cite {trockman2021orthogonalizing,tsuzuku2018lipschitz,prach20241}).}

%% file: acknowledgements.tex
\section*{Acknowledgements}
MJE acknowledges support from the EPSRC (EP/T026693/1, EP/V026259/1, EP/Y037286/1) and the European Union Horizon 2020 research and innovation programme under the Marie Skodowska-Curie grant agreement REMODEL. DM acknowledges support from the EPSRC programme grant in ‘The Mathematics of Deep Learning’, under the project EP/V026259/1. FS acknowledges support from the EPSRC advanced career fellowship EP/V029428/1

%% file: main.bbl
\begin{thebibliography}{10}

\bibitem{antun2020instabilities}
Vegard Antun, Francesco Renna, Clarice Poon, Ben Adcock, and Anders~C Hansen.
\newblock On instabilities of deep learning in image reconstruction and the potential costs of {AI}.
\newblock {\em Proceedings of the National Academy of Sciences}, 117(48):30088--30095, 2020.

\bibitem{pmlr-v70-arjovsky17a}
Mart{\'{\i}}n Arjovsky, Soumith Chintala, and L{\'{e}}on Bottou.
\newblock Wasserstein {G}enerative {A}dversarial {N}etworks.
\newblock In {\em Proceedings of the 34th International Conference on Machine Learning, {ICML} 2017, Sydney, NSW, Australia, 6-11 August 2017}, 2017.

\bibitem{arridge2019solving}
Simon Arridge, Peter Maass, Ozan {\"O}ktem, and Carola-Bibiane Sch{\"o}nlieb.
\newblock Solving inverse problems using data-driven models.
\newblock {\em Acta Numerica}, 28:1--174, 2019.

\bibitem{benning2018modern}
Martin Benning and Martin Burger.
\newblock Modern regularization methods for inverse problems.
\newblock {\em Acta Numerica}, 27:1--111, 2018.

\bibitem{bungert2021clip}
Leon Bungert, Ren{\'e} Raab, Tim Roith, Leo Schwinn, and Daniel Tenbrinck.
\newblock {CLIP}: {C}heap {L}ipschitz {T}raining of {N}eural {N}etworks.
\newblock In {\em International Conference on Scale Space and Variational Methods in Computer Vision}, pages 307--319. Springer, 2021.

\bibitem{celledoniDynamicalSystemsBased2023}
Elena Celledoni, Davide Murari, Brynjulf Owren, Carola-Bibiane Sch{\"o}nlieb, and Ferdia Sherry.
\newblock Dynamical {{Systems}}--{{Based Neural Networks}}.
\newblock {\em SIAM Journal on Scientific Computing}, 45(6):A3071--A3094, 2023.

\bibitem{chen2018neural}
Tian~Qi Chen, Yulia Rubanova, Jesse Bettencourt, and David Duvenaud.
\newblock Neural {O}rdinary {D}ifferential {E}quations.
\newblock In {\em Advances in Neural Information Processing Systems 31: Annual Conference on Neural Information Processing Systems 2018, NeurIPS 2018, December 3-8, 2018, Montr{\'{e}}al, Canada}, pages 6572--6583, 2018.

\bibitem{englRegularizationInverseProblems2000}
Heinz~W. Engl, Martin {Hanke}, and Andreas Neubauer.
\newblock {\em Regularization of {I}nverse {P}roblems}.
\newblock Number 375 in Mathematics and Its Applications. Kluwer Academic Publishers, 2000.

\bibitem{fukushima1980neocognitron}
Kunihiko Fukushima.
\newblock Neocognitron: {A} self-organizing neural network model for a mechanism of pattern recognition unaffected by shift in position.
\newblock {\em Biological Cybernetics}, 36(4):193--202, 1980.

\bibitem{galimberti2023hamiltonian}
Clara~Luc{\'\i}a Galimberti, Luca Furieri, Liang Xu, and Giancarlo Ferrari-Trecate.
\newblock Hamiltonian {D}eep {N}eural {N}etworks {G}uaranteeing {N}onvanishing {G}radients by {D}esign.
\newblock {\em IEEE Transactions on Automatic Control}, 68(5):3155--3162, 2023.

\bibitem{goodfellow2014explaining}
Ian~J. Goodfellow, Jonathon Shlens, and Christian Szegedy.
\newblock Explaining and {H}arnessing {A}dversarial {E}xamples.
\newblock In {\em 3rd International Conference on Learning Representations, {ICLR} 2015, San Diego, CA, USA, May 7-9, 2015, Conference Track Proceedings}, 2015.

\bibitem{gouk2021regularisation}
Henry Gouk, Eibe Frank, Bernhard Pfahringer, and Michael~J Cree.
\newblock Regularisation of {N}eural {N}etworks by {E}nforcing {L}ipschitz {C}ontinuity.
\newblock {\em Machine Learning}, 110:393--416, 2021.

\bibitem{hairer2011geometric}
Ernst Hairer, Christian Lubich, and Gerhard Wanner.
\newblock {\em Geometric {N}umerical {I}ntegration. {S}tructure-{P}reserving {A}lgorithms for {O}rdinary {D}ifferential {E}quations}.
\newblock Springer, 2006.

\bibitem{hairer1993solving}
Ernst Hairer, Syvert~P N{\o}rsett, and Gerhard Wanner.
\newblock {\em Solving {O}rdinary {D}ifferential {E}quations {I}. {N}onstiff {P}roblems}, volume~8.
\newblock Springer, 1993.

\bibitem{hasannasab2020parseval}
Marzieh Hasannasab, Johannes Hertrich, Sebastian Neumayer, Gerlind Plonka, Simon Setzer, and Gabriele Steidl.
\newblock Parseval {P}roximal {N}eural {N}etworks.
\newblock {\em Journal of Fourier Analysis and Applications}, 26:1--31, 2020.

\bibitem{he2015deepresiduallearningimage}
Kaiming He, Xiangyu Zhang, Shaoqing Ren, and Jian Sun.
\newblock Deep {R}esidual {L}earning for {I}mage {R}ecognition.
\newblock In {\em 2016 {IEEE} Conference on Computer Vision and Pattern Recognition, {CVPR} 2016, Las Vegas, NV, USA, June 27-30, 2016}, pages 770--778. {IEEE} Computer Society, 2016.

\bibitem{hertrich2021convolutional}
Johannes Hertrich, Sebastian Neumayer, and Gabriele Steidl.
\newblock Convolutional {P}roximal {N}eural {N}etworks and {P}lug-and-{P}lay {A}lgorithms.
\newblock {\em Linear Algebra and its Applications}, 631:203--234, 2021.

\bibitem{huang2022fi}
Yujia Huang, Ivan Dario~Jimenez Rodriguez, Huan Zhang, Yuanyuan Shi, and Yisong Yue.
\newblock {FI-ODE}: {C}ertifiably {R}obust {F}orward {I}nvariance in {N}eural {ODE}s.
\newblock {\em arXiv preprint arXiv:2210.16940}, 2022.

\bibitem{jaffeLearningNeuralContracting2024}
Sean Jaffe, Alexander Davydov, Deniz Lapsekili, Ambuj~K. Singh, and Francesco Bullo.
\newblock Learning {N}eural {C}ontracting {D}ynamics: {E}xtended {L}inearization and {G}lobal {G}uarantees.
\newblock In {\em Advances in Neural Information Processing Systems 38: Annual Conference on Neural Information Processing Systems 2024, NeurIPS 2024, Vancouver, BC, Canada, December 10 - 15, 2024}, 2024.

\bibitem{kang2021stable}
Qiyu Kang, Yang Song, Qinxu Ding, and Wee~Peng Tay.
\newblock Stable neural {ODE} with {Lyapunov}-stable equilibrium points for defending against adversarial attacks.
\newblock In {\em Advances in Neural Information Processing Systems 34: Annual Conference on Neural Information Processing Systems 2021, NeurIPS 2021, December 6-14, 2021, virtual}, pages 14925--14937, 2021.

\bibitem{kidger2022neural}
Patrick Kidger.
\newblock On {N}eural {D}ifferential {E}quations.
\newblock {\em arXiv preprint arXiv:2202.02435}, 2022.

\bibitem{kolter2019learning}
J.~Zico Kolter and Gaurav Manek.
\newblock Learning {S}table {D}eep {D}ynamics {M}odels.
\newblock In {\em Advances in Neural Information Processing Systems 32: Annual Conference on Neural Information Processing Systems 2019, NeurIPS 2019, December 8-14, 2019, Vancouver, BC, Canada}, pages 11126--11134, 2019.

\bibitem{leimkuhler2004simulating}
Benedict Leimkuhler and Sebastian Reich.
\newblock {\em Simulating {H}amiltonian {D}ynamics}.
\newblock Number~14 in Cambridge Monographs on Applied and Computational Mathematics. Cambridge University Press, 2004.

\bibitem{lunz2018adversarial}
Sebastian Lunz, Carola Sch{\"{o}}nlieb, and Ozan {\"{O}}ktem.
\newblock Adversarial {R}egularizers in {I}nverse {P}roblems.
\newblock In {\em Advances in Neural Information Processing Systems 31: Annual Conference on Neural Information Processing Systems 2018, NeurIPS 2018, December 3-8, 2018, Montr{\'{e}}al, Canada}, pages 8516--8525, 2018.

\bibitem{maass2019deep}
Peter Maass.
\newblock Deep {L}earning for {T}rivial {I}nverse {P}roblems.
\newblock In {\em Compressed Sensing and Its Applications: {{Third}} International {{MATHEON}} Conference 2017}, pages 195--209. Springer International Publishing, Cham, 2019.

\bibitem{madry2017towards}
Aleksander Madry, Aleksandar Makelov, Ludwig Schmidt, Dimitris Tsipras, and Adrian Vladu.
\newblock Towards {D}eep {L}earning {M}odels {R}esistant to {A}dversarial {A}ttacks.
\newblock In {\em 6th International Conference on Learning Representations, {ICLR} 2018, Vancouver, BC, Canada, April 30 - May 3, 2018, Conference Track Proceedings}, 2018.

\bibitem{meunier2022dynamical}
Laurent Meunier, Blaise Delattre, Alexandre Araujo, and Alexandre Allauzen.
\newblock A dynamical system perspective for {Lipschitz} neural networks.
\newblock In {\em International Conference on Machine Learning, {ICML} 2022, 17-23 July 2022, Baltimore, Maryland, {USA}}, pages 15484--15500, 2022.

\bibitem{miyatoSpectralNormalizationGenerative2018}
Takeru Miyato, Toshiki Kataoka, Masanori Koyama, and Yuichi Yoshida.
\newblock Spectral {N}ormalization for {G}enerative {A}dversarial {N}etworks.
\newblock In {\em 6th International Conference on Learning Representations, {ICLR} 2018, Vancouver, BC, Canada, April 30 - May 3, 2018, Conference Track Proceedings}, 2018.

\bibitem{prach20241}
Bernd Prach, Fabio Brau, Giorgio Buttazzo, and Christoph~H Lampert.
\newblock 1-{L}ipschitz {L}ayers {C}ompared: {M}emory, {S}peed, and {C}ertifiable {R}obustness.
\newblock In {\em Proceedings of the IEEE/CVF Conference on Computer Vision and Pattern Recognition}, pages 24574--24583, 2024.

\bibitem{rosenblatt1958perceptron}
Frank Rosenblatt.
\newblock The perceptron: {A} probabilistic model for information storage and organization in the brain.
\newblock {\em Psychological Review}, 65(6):386, 1958.

\bibitem{ryu2019plug}
Ernest~K. Ryu, Jialin Liu, Sicheng Wang, Xiaohan Chen, Zhangyang Wang, and Wotao Yin.
\newblock Plug-and-{P}lay {M}ethods {P}rovably {C}onverge with {P}roperly {T}rained {D}enoisers.
\newblock In {\em Proceedings of the 36th International Conference on Machine Learning, {ICML} 2019, 9-15 June 2019, Long Beach, California, {USA}}, 2019.

\bibitem{Sanz-Serna_1992}
Jesús~María Sanz-Serna.
\newblock Symplectic integrators for {H}amiltonian problems: {A}n overview.
\newblock {\em Acta Numerica}, 1:243–286, 1992.

\bibitem{sherryDesigningStableNeural2024}
Ferdia Sherry, Elena Celledoni, Matthias~J. Ehrhardt, Davide Murari, Brynjulf Owren, and Carola-Bibiane Sch{\"o}nlieb.
\newblock Designing stable neural networks using convex analysis and {{ODEs}}.
\newblock {\em Physica D: Nonlinear Phenomena}, page 134159, 2024.

\bibitem{strogatz2018nonlinear}
Steven~H. Strogatz.
\newblock {\em Nonlinear {{Dynamics}} and {{Chaos}}}.
\newblock CRC Press, 2018.

\bibitem{szegedy2013intriguing}
Christian Szegedy, Wojciech Zaremba, Ilya Sutskever, Joan Bruna, Dumitru Erhan, Ian~J. Goodfellow, and Rob Fergus.
\newblock Intriguing properties of neural networks.
\newblock In {\em 2nd International Conference on Learning Representations, {ICLR} 2014, Banff, AB, Canada, April 14-16, 2014, Conference Track Proceedings}, 2014.

\bibitem{doi:10.1137/1.9781611975161}
Lloyd~N. Trefethen, Ásgeir Birkisson, and Tobin~A. Driscoll.
\newblock {\em Exploring {ODEs}}.
\newblock Society for Industrial and Applied Mathematics, 2017.

\bibitem{trockman2021orthogonalizing}
Asher Trockman and J.~Zico Kolter.
\newblock Orthogonalizing {C}onvolutional {L}ayers with the {Cayley} {T}ransform.
\newblock In {\em 9th International Conference on Learning Representations, {ICLR} 2021, Virtual Event, Austria, May 3-7, 2021}, 2021.

\bibitem{tsuzuku2018lipschitz}
Yusuke Tsuzuku, Issei Sato, and Masashi Sugiyama.
\newblock Lipschitz-margin training: Scalable certification of perturbation invariance for deep neural networks.
\newblock In {\em Advances in Neural Information Processing Systems 31: Annual Conference on Neural Information Processing Systems 2018, NeurIPS 2018, December 3-8, 2018, Montr{\'{e}}al, Canada}, pages 6542--6551, 2018.

\bibitem{yamashita2018convolutional}
Rikiya Yamashita, Mizuho Nishio, Richard Kinh~Gian Do, and Kaori Togashi.
\newblock Convolutional neural networks: an overview and application in radiology.
\newblock {\em Insights into Imaging}, 9:611--629, 2018.

\end{thebibliography}
